\documentclass[11pt]{amsart}
\usepackage{amsmath, amssymb, amscd, amsthm, mathrsfs, url ,pinlabel, verbatim, lipsum, wrapfig}
\usepackage[text={6.5in,9.5in}, centering, letterpaper, dvips]{geometry}
\usepackage{color, multirow, dcpic, latexsym, pictexwd, graphicx, epstopdf, tikz-cd, wrapfig, float, xcolor}
\usepackage{verbatim, relsize, latexsym, latexsym, tikz, mdframed, float, mathtools, flowchart, multicol, stmaryrd, caption, capt-of, old-arrows, dirtytalk, footmisc, setspace, enotez, tabularx, pifont, tensor, csquotes}
\usepackage{svg}
\usepackage{fancyhdr, tikz, tikz-cd, setspace, tensor, enumitem, dirtytalk, nicematrix, dbnsymb}

\usepackage{pst-node}
\usetikzlibrary{positioning, arrows.meta, knots, braids}

\usepackage{scalerel,stackengine}
\stackMath
\newcommand\reallywidehat[1]{%
\savestack{\tmpbox}{\stretchto{%
  \scaleto{%
    \scalerel*[\widthof{\ensuremath{#1}}]{\kern-.6pt\bigwedge\kern-.6pt}%
    {\rule[-\textheight/2]{1ex}{\textheight}}
  }{\textheight}%
}{0.5ex}}%
\stackon[1pt]{#1}{\tmpbox}%
}

\newenvironment{myepigraph}
  {\par\hfill\itshape
   \begin{tabular}{@{}r@{\hspace{0em}}}} 
  {\end{tabular}\par\medskip}

\usepackage[backref=page, linktocpage=true]{hyperref}

\renewcommand*{\backref}[1]{}
\renewcommand*{\backrefalt}[4]{%
    \ifcase #1 (Not cited.)%
    \or        (Cited on page~#2.)%
    \else      (Cited on pages~#2.)%
    \fi}

\definecolor{maroon}{rgb}{0.5, 0.0, 0.0}
\definecolor{darkblue}{rgb}{0.03, 0.27, 0.49}

\hypersetup{
  colorlinks   = true, 
  urlcolor     = maroon, 
  linkcolor    = darkblue, 
  citecolor   = maroon 
}
\usepackage{titletoc}
\usepackage{svg}

\newcommand{\CommaPunct}{\mathpunct{\raisebox{0.5ex}{,}}}

\usepackage[T1]{fontenc}
\usepackage{lmodern}
\DeclareUrlCommand{\bfurl}{}

\newcommand\TB{\overline{tb}}
\newcommand\tb{tb}

\theoremstyle{definition}

\newtheorem{introthm}{Theorem}

\newtheorem{introconj}[introthm]{Conjecture}
\newtheorem{introprob}{Problem}

\theoremstyle{remark}

\newcommand\Z{\mathbb{Z}}
\newcommand\Q{\mathbb{Q}}

\newcommand\C{\mathcal{C}}

\usepackage[colorinlistoftodos,textwidth=3cm]{todonotes}

\newlength\Colsep
\title{Turk's head knots and links: a survey}

\author{Alessio Di Prisa}
\address{Scuola Normale Superiore, 56126 Pisa, Italy}
\email{\url{alessio.diprisa@sns.it}}
\urladdr{\url{https://sites.google.com/view/alessiodiprisa}}

\author{O{\u{g}}uz \c{S}avk}
\address{Department of Mathematics, Middle East Technical University, 06800 \c{C}ankaya, Ankara, Turkey}
\email{\url{savk@metu.edu.tr}}
\urladdr{\url{https://sites.google.com/view/oguzsavk}}

\date{}

\begin{document}

\begin{abstract}

We collect and discuss various results on an important family of knots and links called \emph{Turk's head knots and links} $Th (p,q)$. In the mathematical literature, they also appear under different names such as \emph{rosette knots and links} or \emph{weaving knots and links}. Unless being the unknot or the alternating torus links $T(2,q)$, the Turk's head links $Th (p,q)$ are all known to be alternating, non-split, prime, fibered, hyperbolic, invertible, and periodic. The Turk's head links $Th (p,q)$ are also both positive and negative amphichiral if $p$ is chosen to be odd. Moreover, we highlight and present several more results, focusing on Turk's head knots $Th (3,q)$. We finally list several open problems and conjectures for Turk's head knots and links. We conclude with a short appendix on torus knots and links, which might be of independent interest.
\end{abstract}

\maketitle

\begin{spacing}{1.12}
\setcounter{tocdepth}{1}
\tableofcontents
\end{spacing}

\small{
\begin{myepigraph}
\say{There is no knot with a wider field of usefulness.}\\[1.5ex]
Chapter 17: The Turk's-Head \\
The Ashley Book of Knots\\[3ex]
\end{myepigraph}
}

\vspace{2em}

\section{Objective}

This survey article is addressed to the entire geometry and topology community. The reader should only be familiar with the basics of knot theory; in particular, we follow three standard reference books by Rolfsen \cite{Rol76}, Kawauchi \cite{Kaw90}, and Lickorish \cite{Lic97}. Although the focus on a specific family of knots and links as the subject of a survey may seem limited at first sight, it is nevertheless rich enough to encompass several research articles in the literature, as the extensive bibliography of the paper testifies. We also discuss the other fundamental geometric and topological concepts for Turk's head links and raise various open problems and conjectures to enrich the existing literature, inspired by the recent research articles. We hope that our efforts will have a positive impact and motivate the community to explore their properties in the future.

\section{Introduction}
\label{sec:introduction}

Torus knots and links are one of the most important and central objects in knot theory, low-dimensional topology, and contact geometry. One main reason is Thurston's breakthrough \cite{Thu82, Thu97}, proving that every non-trivial knot in $S^3$ is either a hyperbolic link, or a satellite link, or a torus link. Another crucial reason goes back a century. It is well known that torus links are also realized as the links of Brauner type complex curve singularities $x^p + y^q = 0$, see Milnor's book and references therein \cite[{\sc\S}1]{Mil68}. Since then, several properties of torus links have been well studied, see for instance; the celebrated results of Moser \cite{Mos71}, Litherland \cite{Lit79}, Gordon, Litherland and Murasugi \cite{GLM81}, Kronheimer and Mrowka \cite{KM93, KM95}, Etnyre and Honda \cite{EH01}, and Casals and Gao \cite{CG22}. 

\begin{figure}[htbp]
\centering    
\includegraphics[width=0.55\linewidth]{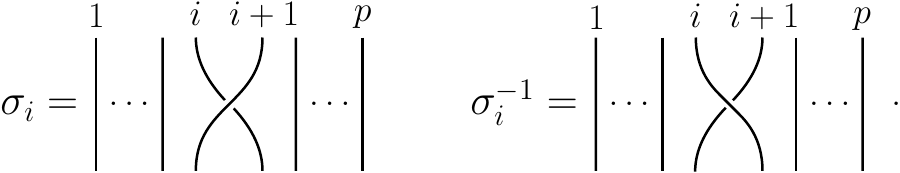}
\caption{The braids $\sigma_i$ and $\sigma^{-1}_i$.}
\label{fig:braids}
\end{figure}

Although torus knots and links have several different but equivalent beautiful descriptions, for our purposes we choose to define them as braid closures. Given two positive integers $p$ and $q$, the \emph{torus link} $T(p,q)$ in $S^3$ is defined as the braid closure of the braid $( \sigma_1 \sigma_2 \sigma_3 \sigma_4  \ldots \sigma_{p-1} )^q$, i.e., $$ T (p,q) \doteq \reallywidehat{ ( \sigma_1 \sigma_2 \sigma_3 \sigma_4  \ldots \sigma_{p-1} )^q } , $$ where the $i$th Artin generator $\sigma_i$ for $i \in \{ 1, \ldots, p-1\}$ of the braid group on $p$ strands $B_p$, and its inverse $\sigma^{-1}_i$ are given in Figure~\ref{fig:braids}.

\begin{figure}[htbp]
\centering

\begin{tikzpicture}
\node[anchor=south west,inner sep=0] at (0,0) {
\includegraphics[width=0.5\columnwidth]{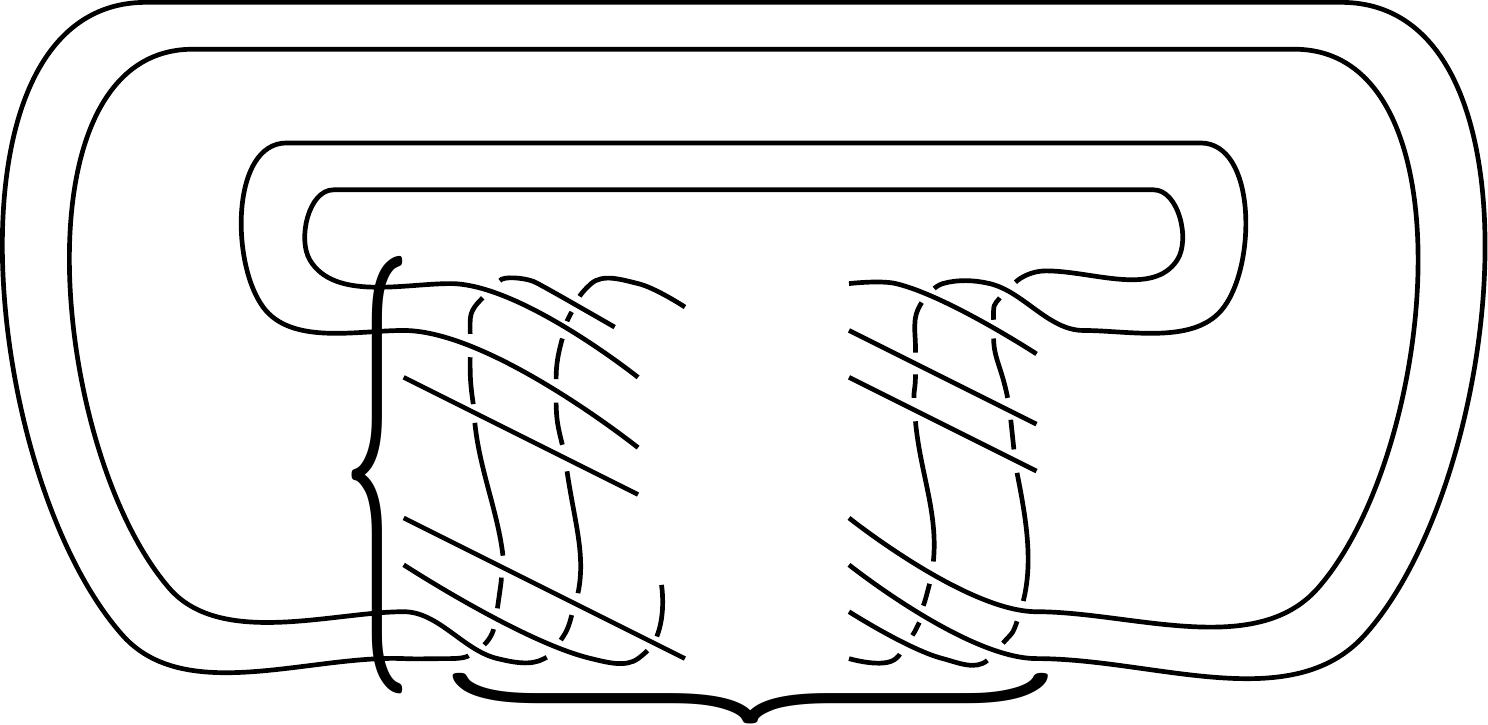}
};
\node[scale=1.5] at (1.7,1.4) (b1){$p$};
\node[scale=1.5] at (4.2,1.4) (b1){$\dots$};
\node[scale=1.5] at (4.2,-0.3) (b1){$q$};
\end{tikzpicture}
\caption{The link diagram of $T(p,q)$.}
\label{fig:torusknots}
\end{figure}

\noindent Using its braid presentation, we can draw a diagram for the torus link $T(p,q)$ as in Figure~\ref{fig:torusknots}. It is well-known that the torus link $T(p,q)$ has $\mathrm{gcd}(p,q)$ components where $\mathrm{gcd}$ denotes the greatest common divisor. When $\mathrm{gcd} (p,q) = 1$, $T(p,q)$ becomes the so-called \emph{$(p,q)$-torus knot}. See {\sc\S}\ref{sec:appendix} for more details.

\begin{figure}[htbp]
\centering

\begin{tikzpicture}
\node[anchor=south west,inner sep=0] at (0,0) {
\includegraphics[width=0.5\columnwidth]{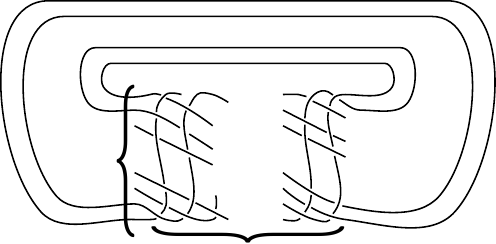}
};
\node[scale=1.5] at (1.7,1.4) (b1){$p$};
\node[scale=1.5] at (4.2,1.4) (b1){$\dots$};
\node[scale=1.5] at (4.2,-0.3) (b1){$q$};
\end{tikzpicture}
\caption{The link diagram of $Th (p,q)$.}
\label{fig:ourknots}
\end{figure}

Let $p$ and $q$ be two positive integers. The $p$-lead $q$-bight \emph{Turk's head link} $Th(p,q)$ in $S^3$ is given by the following braid closure $$ Th (p,q) \doteq \reallywidehat{ ( \sigma_1 \sigma_2^{-1} \sigma_3 \sigma_4^{-1} \ldots \sigma_{p-1}^{\pm 1} )^q }$$ for which
a diagram is clearly obtained by modifying the torus link $T(p,q)$ diagram by replacing the braids with even indices by their inverses, see Figure~\ref{fig:ourknots}. It follows directly from the definition that $Th(p,q)$ is always alternating, for every $p$ and $q$, unlike the corresponding torus link. Similar to the torus link, $Th (p,q)$ has $\mathrm{gcd}(p,q)$ components in total, and it is a knot if $\mathrm{gcd} (p,q) = 1$.

Notice also that as the positive integers $p$ and $q$ vary, we have the following non-trivial intersection $$\{ T(p,q) \} \cap \{ Th (p,q) \} = \{ T(1,q), T(2,q), T(p,1), T(p,2)\} $$ where $T(1,q)$ and $T(p,1)$ are the unknot for any $p \geq 1$ and $q \geq 1$, and $T(2,q)=T(q,2)$ is the only alternating torus link for any $q \geq 2$. When $p \geq 3$ and $q \geq 2$, the Turk's head links $Th (p,q)$ form a very interesting family on their own, which is the subject of this survey. We will see that they are alternating, fibered, hyperbolic, invertible, non-split, periodic, and prime. They are also both positive and negative amphichiral if $p$ is odd.

By switching crossings in their diagrams arbitrarily, one obtains a vast generalization of both torus links and Turk's head links, called \emph{quasitoric links}, which are defined as the braid closures of the \emph{quasitoric braids} of type $(p,q)$ (see also {\sc\S}\ref{sec:rosette}): $$\underbrace{( \sigma_1^{\epsilon_{1,1} } \sigma_2^{\epsilon_{1,2} } \ldots \sigma_{p-1}^{\epsilon_{1,p-1} } ) ( \sigma_1^{\epsilon_{2,1} } \sigma_2^{\epsilon_{2,2} } \cdots \sigma_{p-1}^{\epsilon_{2,p-1} } ) \ldots ( \sigma_1^{\epsilon_{q,1} } \sigma_2^{\epsilon_{q,2} } \ldots \sigma_{p-1}^{\epsilon_{q,p-1} } )}_{q\text{-many}}$$ where $p,q \geq 1$ and $\epsilon_{i,j} = \pm 1$. Manturov and Lamm independently prove that the quasitoric links are universal in the sense that every link in $S^3$ is realized as a closure of a quasitoric braid \cite{Man02, Lam12, Lam99}. 

\subsection*{Organization and Structure} We begin our survey with a detailed literature review in {\sc\S}\ref{sec:literature}. In {\sc\S}\ref{sec:general}, we discuss the Turk's head knots and links $Th (p,q)$ in the general case. We describe their rich properties and present various results on their invariants. In {\sc\S}\ref{sec:special}, we focus on the Turk's head knots $Th (3,q)$, assuming $\mathrm{gcd}(3,q) = 1$. We want to emphasize that the literature in this special case is extensive, we therefore share several additional results. We list several open problems related to the Turk's head knots and links in {\sc\S}\ref{sec:problems}. Finally, {\sc\S}\ref{sec:appendix} is an appendix on the torus knots and links. 

We make several subsections to engage the readers within the survey, whose names are clearly in line with the properties and invariants of our knots and links. Moreover, we create two additional table of contents for {\sc\S}\ref{sec:general} and {\sc\S}\ref{sec:special} to specify their subsections.


\section{Literature Review}
\label{sec:literature}

In this section, we list and review various articles about the Turk's head knots and links $Th (p,q)$. The order of the following subsections also reflects their chronology in the literature.

\subsection{Turk's Head Knots and Links}

The earliest record of Turk's head links in the mathematical literature is Fox's 1961 article \cite{Fox62}, in which Fox posed a problem about the periodic symmetries of knots, which he called the \emph{Turk's head problem}. This problem was later partially solved by Murasugi \cite[Corollary~1]{Mur71}. In his paper, Fox specifically emphasized that such knots and links are traditionally called \emph{Turk's head knots and links} since they appear in the literature earlier, see \cite[p.~182]{Fox62}. 

Conway later reintroduced Turk's head links using his new notation $((p-1) \times q)^* $, \cite[p.~341]{Con70}. Conway and Trotter were interested in Turk's head knots from the viewpoint of sliceness. Trotter in particular observed that several concordance invariants vanish simultaneously for specific Turk's head knots and do not give obstructions. Focusing on Turk's head knots, Conway also addressed the slice-ribbon conjecture.

Certain examples of Turk's head links also appear in the books by Crowell and Fox \cite[p.~90 and p.~126]{CF63} and Rolfsen \cite[p.~338]{Rol76}, the unfinished book project by Bonahon and Siebenmann \cite{BS79}, Thistlethwaite's survey on knot tabulations \cite{Thi85}, and in Kirby's 1997 problem list \cite[Problem~1.62]{Kir78}. 

A widely known classical reference for Turk's head knots and links is Ashley's famous encyclopedia of knots, which contains a comprehensive discussion, see \cite[{\sc\S}17]{Ash44}. Various illustrative descriptions of Turk's head links in Ashley's book seem to have been popularized among topologists with the help of the book by Brieskorn and Kn\"orrer, which contains some scanned pages from Ashley's encyclopedia, see \cite[pp.~449-454]{BK86}. A more recent historical reference for a thorough discussion of Turk's head knots and links is the book edited by Turner and van de Griend \cite{TG96}. In addition to the several enlightening illustrations, the book contains an anthology of various essays with cultural and historical perspectives.

The name \say{Turk's head} derives from the old Turkish calpacs, i.e., the traditional high-crowned caps worn by Turkish men during the times of the Ottoman Empire. The earliest written references in the literature date back to the late eighteenth century; see the recent survey article by A. {\AA}str{\"o}m and C. {\AA}str{\"o}m \cite{AA21}, and the references therein, for a further historical discussion with detailed archaeological aspects. See also Pennock's work \cite{Pen05} for their various excellent illustrations.

By the turn of the millennium, Turk's head links appeared in many research articles and have been studied from several different perspectives; see for instance, the articles by Bozh{\"u}y{\"u}k \cite{Boz80, Boz80b, Boz85, Boz93}, Harer \cite{Har82}, Edmonds and Livingston \cite{EL83}, Long \cite{Lon83, Lon86}, Kojima and Long \cite{KL88}, Meyerhoff and Ruberman \cite{MR90}, Mednykh and Vesnin \cite{MV95a, MV95b, MV96}, Jones and Przytycki \cite{JP98}, and Hilden, Lozano and Montesinos \cite{HLM00}. They were also studied in the dissertations of Long \cite{Lon83} and Bergbauer \cite{Ber98}. 

The work of Nakanishi and Yamada in 2000 \cite{NY00} made the study of Turk's head links attractive once again among topologists. The basic properties of Turk's head links $Th (p,q)$ were studied and reviewed in \cite{NY00}. Since then, their properties have been further studied in the articles by Dowdall, Mattman, Meek, and Solis \cite{DM10}, Lopes and Matias \cite{LM15}, and Takemura \cite{Tak16, Tak18}. All these results are presented in detail in the following sections.

\subsection{Rosette Knots and Links}
\label{sec:rosette}

Three years after Fox's article \cite{Fox62}, Kr{\"o}tenheerdt introduced the Turk's head links as the \emph{rosette links} (\emph{Rosettenknoten} in German) in his first article, which was written in German \cite{Kro64}. The links were therein defined by using braid closures, as we do here. The name \say{rosette} is a reference to the fact that Turk's head links look like ornaments in the shape of roses.

The main purpose of \cite{Kro64} was to resolve the isotopy classification problem of the links $Th (p,q)$ in terms of $p$ and $q$. Along the way, Kr{\"o}tenheerdt observed that $Th (p,q)$ is amphichiral (i.e. it is isotopic to its mirror image) if $p$ is odd, and provided a partial solution to this problem. 

Later, the complete classification result was given by Murasugi \cite{Mur65} using the Trotter--Murasugi signatures. The key observation was that $Th (p,q)$ can be written as the Murasugi sum of copies of certain torus links. A few years later, Kr{\"o}tenheerdt reproved this result by analyzing the coefficients of Alexander polynomials of $Th (p,q)$, see \cite{Kro71}. We explain the details of these results below.

The most recent references using the notion of rosette links are the articles by Lamm and Obermeyer \cite{LO99}, Lamm \cite{Lam12, Lam99}, and Koseleff and Pecker \cite{KP14}. The former showed that these links are all billiard knots in a cylinder, i.e. they are periodic billiard trajectories without self-intersections in some cylindrical billiard room. In \cite{Lam12, Lam99}, Lamm used Turk's head links in the sense of quasitoric links and proved the universality result we presented in the introduction. Koseleff and Pecker proved a conjecture of Jones and Przytycki \cite{JP98} by showing that every link is a billiard link in a convex right prism, see \cite{KP14}. They also noted that their proof follows from \cite{LO99} for knots.

\subsection{Weaving Knots and Links}

More recently, Turk's head links were introduced as the \emph{weaving links} in the influential articles by Champanerkar, Kofman, and Purcell \cite{CKP16a, CKP16b}. They were extensively studied as hyperbolic links, which X.-S Lin \cite{Lin05}previously suggested. It turns out that Turk's head links are diagrammatically and geometrically maximal, in the sense that their determinants and volumes are asymptotically maximal with respect to their crossing numbers.

\begin{figure}[htbp]
    \centering    
    \includegraphics[width=0.15\linewidth]{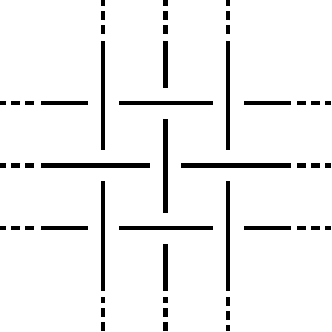}
    \caption{The infinite weave $\mathcal{W}$.}
    \label{fig:weave}
\end{figure}

When Champanerkar, Kofman, and Purcell presented weaving links, their similarity to the torus links was highlighted. A motivation for the name \say{weaving} is that these links are related to the \emph{infinite weave} $\mathcal{W}$ in Figure~\ref{fig:weave}, which is the infinite alternating link with the square grid projection. The complement $\mathbb{R}^3 \setminus \mathcal{W}$ has a complete hyperbolic structure obtained by tessellating the manifold by regular ideal octahedra. 

Since then, Turk's head links have been studied under the name of weaving links by various topologists; see, for instance, the articles by Burton \cite{Bur16}, Jablan, Kauffman, Sazdanovic, Sto{\v s}i{\'c}, and Zekovi{\'c} \cite{JKSSZ16}, Adams, Kaplan-Kelly, Moore, Shapiro, Sridhar, and Wakefield \cite{AKMSSW18}, Mishra and Raundal \cite{MR19}, Owad \cite{Owa19}, Gan \cite{Gan21}, Mishra and Staffeldt \cite{MS21}, Singh, Mishra, and Ramadevi \cite{SMR21}, Champanerkar, Kofman, and Purcell \cite{CKP22}, Singh and Chbili \cite{SC22}, Joshi, Negi and Prabhakar \cite{JNP23}, AlSukaiti and Chbili \cite{AC23}, Kaiser and Mishra \cite{KM24}, Gill, Joshi and Shimizu \cite{GJS24}. In the upcoming sections, we thoroughly discuss these results.


\section{Results on the General Case}
\label{sec:general}

\startcontents[sections]
\setcounter{tocdepth}{2}
\printcontents[sections]{}{1}{}

\vspace{1.5em}

Let $p$ and $q$ be two positive integers. The Turk's head link $Th (p,q)$ has $\mathrm{gcd}(p,q)$ components in general, see for instance the note of Turner and Schaake \cite{TS91}. If $\mathrm{gcd} (p,q) = 1$, then $Th(p,q)$ is a knot.

\subsection{Isotopy Classification}  

First, we describe the classification of Turk's head links up to isotopy. The first partial attempt to classify Turk's head links was due to Kr{\"o}tenheerdt \cite{Kro64}. The classification was completed by Murasugi \cite{Mur65} via Trotter--Murasugi signatures of links. Later, Kr{\"o}tenheerdt \cite{Kro71} and Nakanishi and Yamada \cite{NY00} reproved this result by computing coefficients of Alexander polynomials and Seifert genera of Turk's head links, respectively.

\begin{introthm} We have that
    \begin{itemize}[leftmargin=2em]
    \item $Th (1,q)$ and $Th (p,1)$ are the unknot for any $p$ and $q$,
    \item $Th (2,q)$ is the alternating torus link $T (2,q)$ for any $q$,
    \item $Th (p,q)$ and $Th (p',q')$ have distinct knot types if $p \geq 2$, $q \geq 2$ and $(p,q) \neq (p',q')$.
\end{itemize}
\end{introthm}

\noindent From now on, we assume that $p \geq 3$ and $q \geq 2$.

\subsection{Diagrammatic Properties}
\label{sec:crossing}

By the definition given in {\sc\S}\ref{sec:introduction}, the link $Th (p,q)$ is alternating with $(p-1)q$ crossings in total. In particular, in the link diagram in Figure~\ref{fig:ourknots}, we have that $$ \# \{ \text{positive crossings} \ \overcrossing  \} = \begin{cases} 
      \frac{(p-1)q}{2}, & \text{if} \ p \ \text{is odd,} \\
      \frac{pq}{2}, & \text{if} \ p \ \text{is even,}
       
\end{cases} \ \ \ \ \ \# \{ \text{negative crossings} \ \undercrossing  \} = \begin{cases} 
      \frac{(p-1)q}{2}, & \text{if} \ p \ \text{is odd,} \\
      \frac{(p-2)q}{2}, & \text{if} \ p \ \text{is even.}
\end{cases} $$ Therefore, the writhe of this diagram is given by $$w (Th(p,q)) = \begin{cases} 
      0, & \text{if} \ p \ \text{is odd,} \\
      q, & \text{if} \ p \ \text{is even.}
\end{cases}$$

Due to the resolution of Tait's conjecture independently by Kauffman \cite{Kau87}, Murasugi \cite{Mur87} and Thistlethwaite \cite{Thi87}, the crossing number (denoted by $c$) in the link diagram in Figure~\ref{fig:ourknots} is minimal, hence $c (Th (p,q)) = (p-1)q$.

A link $L \subset S^3$ is called a \emph{split link} if there is a $2$-sphere in $S^3 \setminus L$ separating $S^3$ into two balls, each of which contains at least one component of $L$, see \cite[Definition~4.1]{Lic97}. Otherwise, $L$ is said to be \emph{non-split}. Recall that a diagram realizing the crossing number of a split link $L$ is also split.

Another application of the resolution of Tait's conjecture implies that $Th (p,q)$ is non-split. Alternatively, the same result follows from Kr{\"o}tenheerdt's computation in \cite{Kro71} which yields that the Alexander polynomial of $\Delta_{Th (p,q)} (t)$ is non-zero. This implies that the link $Th (p,q) $ is non-split, see for instance \cite[Proposition~6.3.4]{Mur96}. This fact was recently reproved by Gill, Joshi and Shimizu \cite[Proposition~2.3]{GJS24} by using Menasco's work \cite{Men84}. They also showed that the Turk's head link is \emph{prime}, i.e. $Th (p,q)$ cannot be represented as a connected sum of other links.

The $3$-sphere $S^3$ has a natural open book decomposition whose pages are open disks and binding is an unknotted circle. Then the binding circle is regarded as the standard $z$-axis and each page is considered as a half-plane $\{ (r, \theta, z) \ \vert \ r \geq 0 \}$ in the cylindrical polar coordinate system. Every link $L \subset S^3$ can be embedded in an open book with finitely many pages so that it meets each page in a simple arc. So, such an embedding is called an \emph{arc presentation} of $L$. The \emph{arc index} $\alpha (L)$ of $L$ is the minimum number of pages in any arc presentation for $L$, see the articles by Birman and Menasco \cite{BM94} and Cromwell \cite{Cro95}. The works of Morton and Beltrami \cite{MB98} and Bae and Park \cite{BP00} imply that we have $$\alpha ( Th(p,q) ) = c ( Th(p,q) ) +2 = (p-1)q +2 .$$

Given a diagram for any oriented knot, we can define its \emph{Gauss code} as follows. First, pick arbitrarily a point on the knot away from its crossings. Then, follow the orientation of the knot to arrive at the first crossing and label it 1. Next, follow the strand to the next crossing. If the crossing you arrive at has not already labeled, then label it 2. If not, skip this crossing and repeat the same procedure. Continue until all crossings have been labeled once. Finally, leave the number the same if the crossing is an overcrossing, and change it with a negative sign if it is an undercrossing.

Following the article by Radovi{\'c} and Jablan \cite{RJ15}, an alternating knot is said to be \emph{Gauss code ordered} if it has at least one minimal diagram with an ordered Gauss code. Since every Gauss code ordered knot diagram is completely determined by the second half of its ordered Gauss code \cite{RJ15}, it is said to be a \emph{short Gauss code}. Inspired by the work of Radovi{\'c} and Jablan, Owad introduced the notion of a \emph{straight number} of a knot $K$, $\mathrm{str} (K),$ as the minimum $n$ such that $K$ has an $n$ crossing \emph{straight diagram} \cite{Owa18},\footnote{See \cite[Figure~1]{Owa18} for the straight diagram of the figure-eight knot $Th(3,2)$.} and studied straight numbers of Turk's head knots. Here, given $\sigma \in S_n$ (the symmetric group on $n$ elements), the straight diagram means that a knot diagram has a short Gauss code of the following form $$(1,2,3, \ldots, n-1, n, \sigma(1), \sigma(2), \sigma(3), \ldots , \sigma(n-1), \sigma(n) ).$$ 

Another motivation for \cite{Owa18} was the work of Adams, Shinjo and Tanaka \cite{AST11} in which they proved that every knot admits a diagram with two arcs where all crossings occur between these two arcs. Assuming that one of these arcs is straight in a horizontal way, we say that the diagram is in the straight position. So, the other definition for the \emph{straight number} of a knot is the minimum number of crossings over all diagrams of the knot in straight position. A knot is called \emph{perfectly straight} if $\mathrm{str} (K) = c(K)$. 

We have the following results proven by Owad \cite{Owa18, Owa19}: \begin{itemize}
    \item $Th (3,2)$ is perfectly straight which corresponds to the figure eight knot,
    \item $Th (p,q)$ is not perfectly straight for $p \geq 3$ and $q \geq p+1$. In particular, $\mathrm{str} (Th (p,q)) > c(Th (p,q))$.
\end{itemize}

\subsection{Murasugi Sums and Trotter--Murasugi Signatures}
\label{sec:signatures} 

The \emph{Trotter--Murasugi signature} of a link $L \subset S^3$ is denoted by $\sigma (L)$ and given by the signature of the square matrix $V + V^T$ where $V$ is a Seifert matrix for $L$, see \cite{Tro62, Mur65b}.

\begin{figure}[htbp]
\centering
\includegraphics[width=0.75\columnwidth]{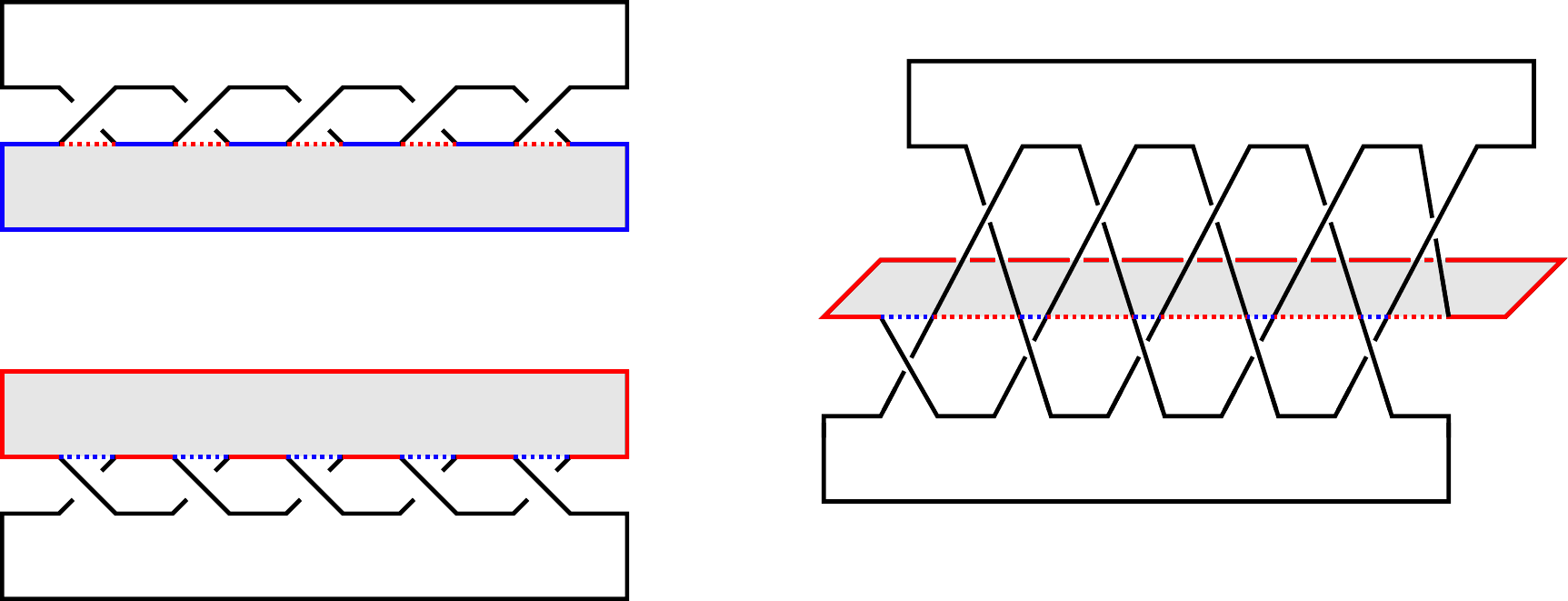}
\caption{The Murasugi sum decomposition: $Th(3,5) = T(2,5) \ \#_{10} \ \overline{T(2,5)}$. The red and blue edges together constitute the gray-shaded $10$-gon inside the Seifert surface. The two $10$-gons are glued by identifying blue edges with blue edges and red edges with red edges.}
\label{fig:murasugisum}
\end{figure}

Now we present an important operation so-called the \emph{Murasugi sum} for the construction of surfaces \cite{Mur58, Mur60}. Given two links $L_1$ and $L_2$ with Seifert surfaces $F_1$ and $F_2$, their Murasugi sum $L_1 \#_{2n} L_2$ is formed by gluing $F_1$ and $F_2$ along an embedded $2n$-gon so that the boundary of the resulting surface $F$ is $L_1 \#_{2n} L_2$, see Figure~\ref{fig:murasugisum}. More precisely, we have
\begin{enumerate}
    \item $F = F_1 \cup F_2$ with $F_1 \cap F_2 = D^2$ such that \begin{itemize}
        \item $\partial D^2$ is a $2n$-gon with edges $a_1, b_1, \ldots, a_n,b_n$,
        \item $a_i$ is contained $\partial F_1$ and is a proper arc in $F_2$ for all $1\leq i \leq n$,
        \item $b_j$ is contained $\partial F_2$ and is a proper arc in $F_1$ for all $1\leq j \leq n$.
    \end{itemize}
   \item[(2)] There exist $3$-balls $B^3_1$ and $B^3_2$ in $S^3$ such that \begin{itemize}
       \item $B^3_1 \cup B^3_2 = S^3$ and $B^3_1 \cap B^3_2 = S^2$,
       \item $F_1 \subset B^3_1$ and $F_2 \subset B^3_2$,
       \item $\partial B^3_1 \cap F_1 = D^2 = \partial B^3_2 \cap F_2$.
   \end{itemize}
\end{enumerate}

In line with Gabai's slogan in the title of his papers \cite{Gab83, Gab85}, we say that the Murasugi sum is a natural operation, and it generalizes the connected sum operation in the sense that $\#_2 = \#$. See \cite[{\sc\S}4.2]{Kaw90} for more details.

\vspace{0.5em}

In \cite{Mur65} (see also \cite[{\sc\S}3]{Mur60}), Murasugi observed that $Th (p,q)$ can be decomposed as the Murasugi sum of $(p-1)$ copies of positive and negative torus links (see {\sc\S}\ref{sec:appendix}). Recall that the Turk's head link is constructed by using positive (resp. negative) bights of the braids $\sigma_i$ (resp. $\sigma^{-1}_i$) for $1 \leq i \leq p-1$ and repeating this process $q$ times consecutively, so we have the following splitting.

\begin{introthm} 
\label{thm:murasugi}
We have $$Th(p,q) = \begin{cases} 
      \underbrace{T(2,q) \ \#_{2q} \ \overline{T(2,q)} \ \#_{2q} \ T(2,q) \ \#_{2q} \ \overline{T(2,q)} \ \#_{2q} \ \ldots \ \#_{2q} \  T(2,q) }_{(p-1)\text{copies}} \CommaPunct & \text{if} \ p \ \text{is even,} \vspace{1 em} \\ 
      \underbrace{T(2,q) \ \#_{2q} \ \overline{T(2,q)} \ \#_{2q} \ T(2,q) \ \#_{2q} \ \overline{T(2,q)} \ \#_{2q} \ \ldots \ \#_{2q} \ \overline{T(2,q)} }_{(p-1)\text{copies}} \CommaPunct & \text{if} \ p \ \text{is odd.}  
\end{cases} $$
\end{introthm}

\vspace{0.5em}

\noindent We call the splitting in Theorem~\ref{thm:murasugi} the \emph{Murasugi sum decomposition}. See Figure~\ref{fig:murasugisum} for a sample construction. Consult also \cite[{\sc\S}3]{Ber98} for a detailed discussion.

The above Murasugi sum decomposition of Turk's head link is very crucial and has profound consequences in this survey. The first important application is the computation of the Trotter--Murasugi signature of the link $Th (p,q)$, see \cite{Mur65} for the first proof. 

\begin{introthm}\label{thm:signatures}
We have $$ \sigma (Th (p,q)) = \sum_{i=1}^{p-1} (-1)^{i+1} (q-1) = \begin{cases} 
      q-1, & \text{if} \ p \ \text{is even,}\\
      0, & \text{if} \ p \ \text{is odd.}  
\end{cases}$$
\end{introthm}

Here, the first equality follows from the additivity of the Trotter--Murasugi signature under the Murasugi sum for alternating links \cite[Theorem~5.4]{Mur65b}. Moreover, $ \sigma (T (2,q) ) = -(q-1)$ and $\sigma (\overline{T (2,q)} ) = q-1 ,$ see {\sc\S}\ref{sec:appendix} for the signature formulas.

\subsection{Symmetries} 
\label{sec:symmetries}

\begin{figure}[htbp]
\centering
\includegraphics[width=0.65\columnwidth]{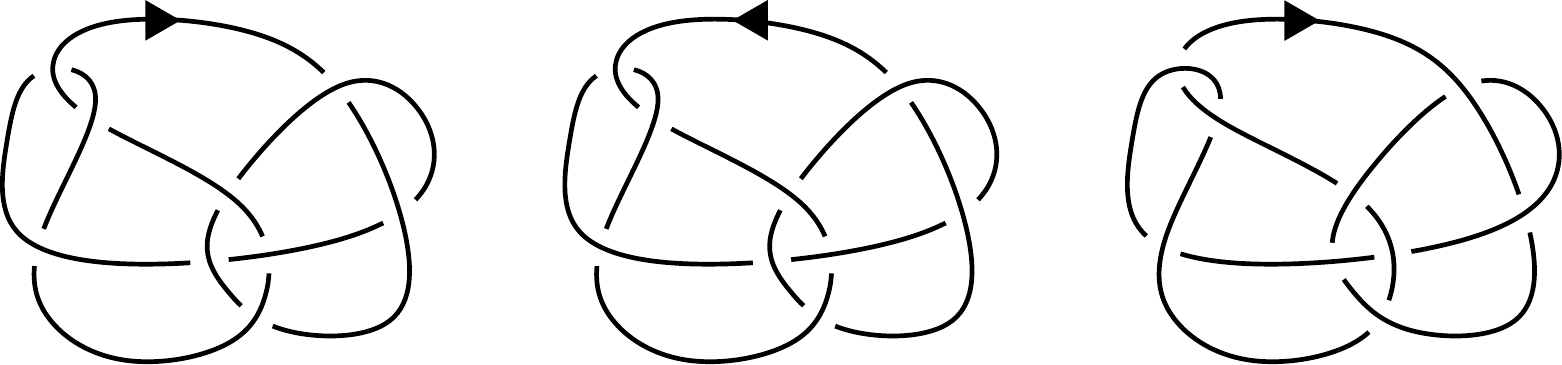}
\caption{From left to right: the knot $K = 9_{32}$, its reverse $-K$ and its mirror image $\overline{K}$. These three knots are different up to isotopy, see KnotInfo \cite{knotinfo}.}
\label{fig:mirrorreverse}
\end{figure}

Given an oriented knot $K$, its \emph{reverse} $-K$ (resp. its \emph{mirror image} $\overline{K}$) is obtained by changing the orientation (resp. the crossings) of $K$, see Figure~\ref{fig:mirrorreverse}.

There are four main symmetries for knots, and the following definitions are naturally extended to links by sending each knot to itself under isotopy maps. We call a knot $K$ in $S^3$: \begin{itemize}
    \item \emph{invertible} if $K$ is isotopic to $-K$. If the isotopy of $S^3$ is given by an involution, $K$ is said to be \emph{strongly invertible}.
    \item \emph{positive amphichiral} if $K$ is isotopic to $\overline{K}$. If the isotopy of $S^3$ is given by a finite order map, $K$ is said to be \emph{periodically positive amphichiral}. If this map can be chosen to be an involution, $K$ is said to be \emph{strongly positive amphichiral}.
    \item \emph{negative amphichiral} if $K$ is isotopic to $-\overline{K}$. If the isotopy of $S^3$ is given by an involution, $K$ is said to be \emph{strongly negative amphichiral}.
    \item (\emph{cyclically})\footnote{There is another notion of periodicity: a knot is said to be \emph{freely $n$-periodic} if it is invariant under a fixed point free map of order $n$. We have omitted it from the discussion since it does not appear elsewhere in the survey.} \emph{$n$-periodic} if $K$ is invariant under a rotation of $2\pi/n$ radians.
\end{itemize}

In \cite{Kro64}, Kr{\"o}tenheerdt proved that $Th (p,q)$ is both positive and negative amphichiral if $p$ is odd. This result was reproved in \cite{NY00}. Nakanishi and Yamada \cite{NY00} also showed that $Th (p,q))$ is always invertible. It is clear from its definition that the Turk's head link $Th (p,q)$ is $q$-periodic. Finally, in \cite[{\sc\S}7]{Lon83} (see also \cite{Lon84}), Long observed that $Th (p,q))$ is strongly positive amphichiral if $p$ and $q$ are both odd.

Since $Th (p,q)$ is always hyperbolic (see {\sc\S}\ref{sec:hyperbolic-general}), Kawauchi's work \cite{Kaw79} implies that $Th (p,q)$ is strongly invertible and strongly negative amphichiral. Moreover, $Th (p,q)$ is periodically positive amphichiral, i.e., the isotopy is given by a finite order map (so an order $2$ isotopy corresponds to an involution). See \cite[{\sc\S}2.1]{DPS24} for more details.

Three important symmetries for the Turk's head knot $Th (3,5)$ are depicted in Figure~\ref{fig:Th(3,5)}. In general, the Turk's head link has precisely the following symmetries in a nutshell.

\begin{introthm}
We have that \begin{itemize}[leftmargin=2em]
    \item $Th (p,q) $ is always $q$-periodic,
    \item $Th (p,q)$ is always strongly invertible, 
    \item $Th (p,q) $ is strongly negative amphichiral if $p$ is odd,
    \item $Th (p,q) $ is periodically positive amphichiral if $p$ is odd,
    \item $Th (p,q) $ is strongly negative and strongly positive amphichiral if $p$ and $q$ are both odd.
\end{itemize}
\end{introthm}

\begin{figure}[htbp]
\centering
\includegraphics[width=0.55\columnwidth]{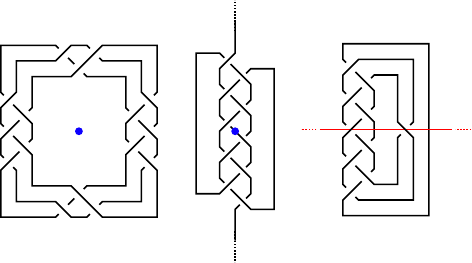}
\caption{From left to right: strongly positive amphichiral, strongly negative amphichiral and strongly invertible. In the first two cases, the involution is given by the $\pi$-rotation around the blue dot composed with the reflection along the plane of the diagram. The third symmetry is given by the $\pi$-rotation around the red axis.}
\label{fig:Th(3,5)}
\end{figure} 

\subsection{Seifert Genera and Fiberedness}
\label{sec:genus-fiberedness}

Given a link $L \subset S^3$, the \emph{Seifert genus} or \emph{$3$-genus} of $L$ is denoted by $g_3 (L)$ and it is given by the minimum genus of any Seifert surface for $L$, see Seifert's article \cite{Sei35}.

A Seifert surface associated with $Th (p,q)$ is depicted in Figure~\ref{fig:murasugisum} for the knot $Th(3,5)$. Nakanishi and Yamada \cite{NY00} indicated that we are able to determine the Seifert genus $g_3$ of $Th (p,q)$ completely by either using Murasugi's result \cite{Mur63} (see also Crowell \cite{Cro59}) or Gabai's result \cite{Gab83}. In particular, we have the following result.

\begin{introthm}
We have $$g_3 (Th (p,q)) = \frac{(p-1)(q-1) + 1 - \mathrm{gcd} (p,q) }{2} \cdot $$ 
\end{introthm}

Another important notion for a link is fiberedness. A link $L \subset S^3$ is called \emph{fibered} if the link complement $S^3 \setminus L$ admits a fibration over the circle $S^1$. As emphasized by Long \cite[{\sc\S}7]{Lon83}, the following result can be obtained by applying Stallings' result \cite[Theorem~1]{Sta78} and using on the Murasugi sum decomposition in Theorem~\ref{thm:murasugi}. Meyerhoff and Ruberman \cite[p.~758]{MR90} also observed the same result by relying on Gabai's work \cite[Theorem~3]{Gab83}. 

\begin{introthm}
For any $p$ and $q$, the Turk's head link $Th (p,q)$ is fibered.
\end{introthm}

The above observations further imply that the Seifert surface obtained by the Murasugi sum description of $Th (p,q)$ (as shown in Figure~\ref{fig:murasugisum}) is minimal. Finally, Kobayashi's work \cite{Kob89} further guarantees that the link $Th (p,q)$ bounds a unique minimal genus Seifert surface relative to the boundary up to isotopy in $S^3$, cf. Kim, Miller and Yoo \cite{KMY24}.

\subsection{Seifert Matrices and Alexander Polynomials}
\label{sec:alexander-general}

We can also describe the Seifert matrix $V_{p,q}$ associated with the unique minimal Seifert surface for $Th (p,q)$ by relying on the results by Kr{\"o}tenheerdt \cite{Kro71} and Takemura \cite{Tak18}. Note that all the empty entries (resp. empty blocks) are zero entries (resp. zero blocks) for the following matrices, and the notation $T$ stands for the transpose operation. 

As another application of the Murasugi sum decomposition in Theorem~\ref{thm:murasugi}, the Seifert matrix $V_{p,q}$ associated to Turk's head link $Th (p,q)$ is given by the following square matrix of size $(p-1)(q-1)$: $$V_{p,q} \doteq \begin{cases} 

\begin{pNiceMatrix}
 A_q    &            &        &       &        &        \\
-A_q   & -A^T_q     &        &        &        &        \\
       & -A_q       & A_q    &        &        &        \\
       &            & \ddots & \ddots &        &        \\
       &            &        & -A_q   & -A^T_q &        \\
       &            &        &        & -A_q   & A_q    \\
\end{pNiceMatrix}, & \ \text{if} \ p \ \text{is even,} \\

\begin{pNiceMatrix}
 A_q    &            &        &       &        &        \\
-A_q   & -A^T_q     &        &        &        &        \\
       & -A_q       & A_q    &        &        &        \\
       &            & \ddots & \ddots &        &        \\
       &            &        & -A_q   & A_q &        \\
       &            &        &        & -A_q   & -A^T_q    \\
\end{pNiceMatrix}, & \ \text{if} \ p \ \text{is odd,} \end{cases}$$ 

\noindent where $A_q$ is the Seifert matrix for the torus link $T (2,q)$, which is the square matrix of size $q-1$: 

$$A_q = \begin{pNiceMatrix}
-1     & 1      &        &        &        \\
       & -1     & 1      &        &        \\
       &        & \ddots & \ddots &        \\
       &        &        & \ddots & 1      \\
       &        &        &        & -1     \\       
\end{pNiceMatrix}.$$ 

\vspace{1em}

Since our link $Th (p,q)$ is alternating and fibered, we make the following two observations due to the classical results of Murasugi \cite{Mur58}, Stallings \cite{Sta61} and Neuwirth \cite{Neu61}: \begin{enumerate}
    \item The Alexander polynomial $\Delta_{Th (p,q)} (t)$ is monic,
    \item The degree of $\Delta_{Th (p,q)} (t)$ is $(p-1)(q-1) + 1 - \mathrm{gcd} (p,q)$, i.e., twice $g_3 (Th (p,q))$.
\end{enumerate}

Now, we describe the Alexander polynomial of $Th (p,q)$, following the article of Takemura \cite{Tak18}. For $q \geq 2$, let $B_q$ be the square matrix of size $q-1$ defined by 

$$B_2 = (1) \ \ \text{and} \ \  B_q = \begin{pNiceMatrix}
0      & \ldots & 0      & 1        \\
-1     & 1      &        & \vdots   \\
       & \ddots &        & 1        \\
       &        & -1     & 1        \\       
\end{pNiceMatrix} \ \text{for} \ q \geq 3.$$ 

\vspace{0.5 em}

\noindent Also, let $X_q (t)$ and $Y_q (t)$ be the matrices defined by the identities $$X_q (t) = t^{\frac{1}{2}} I_{q-1} - t^{\frac{-1}{2}} B_q \ \ \text{and} \ \ Y_q (t) = X_q (t^{-1}),$$ where $I_{q-1}$ is the identity matrix of size $q-1$. Then we define the square matrix $P_{p,q} (t)$ of size $q-1$ inductively as follows: \begin{itemize}
    \item for $p=2,3$, we have $$P_{2,q} (t) = X_q (t) \ \ \text{and} \ \ P_{3,q} (t) = t^{\frac{-1}{2}} X_q (t) Y_q (t) + t^{\frac{-1}{2}} B_q ,$$
    \item for $p \geq 4$, we have $$P_{p,q} (t) = \begin{cases} 
     t^{\frac{-1}{2}} X_q (t) P_{p-1,q} (t) + t^{-1} B_q P_{p-2,q} (t) , & \text{if} \ p \ \text{is even,}\\
     t^{\frac{-1}{2}} Y_q (t) P_{p-1,q} (t) + t^{-1} B_q P_{p-2,q} (t) , & \text{if} \ p \ \text{is odd,}  
\end{cases}$$
\end{itemize}

\noindent so that the Alexander polynomial of our link is given by $$ \Delta_{Th (p,q)} (t) = \mathrm{det} (P_{p,q} (t)). $$ A more refined but still a complicated formula for $P_{p,q} (t)$ can be found in \cite[Lemma~3.2]{Tak18}. 

Now, we state the main result of Takemura's article \cite{Tak18} (cf. \cite{Mur71, Sak81}) about the divisibility of Alexander polynomials, which we call the \emph{Takemura's divisibility criterion}.

\begin{introthm}
\label{thm:takemura}
For any positive integers $a, b, p$ and $q$, we have $$\Delta_{Th (p,q)} (t) \ \big\vert \ \Delta_{Th (ap,bq)} (t) .$$
\end{introthm}

\subsection{Determinants}
\label{sec:determinants}

A knot diagram $D$ is said to be \emph{alternating} if every crossing alternates between under- and over-passes consecutively or vice versa. Moreover, a crossing in $D$ is called \emph{nugatory} if the knot type remains the same after the crossing change on this one.

Next, we present the work of Dowdall, Mattman, Meek, and Solis \cite{DM10} on the Kauffman--Harary conjecture, which expects that given an alternating knot diagram $D$ with no nugatory crossings if the determinant of $D$ is a prime number $p$, then every non-trivial Fox $p$-coloring of $D$ assigns different colors to different arcs of $D$, see \cite{HK99, DM10}. 

They proved the Kauffman--Harary conjecture for all $Th (p,q)$ except for the case where $p \geq 5$ is odd and $q \geq 3$. For this purpose, Dowdall, Mattman, Meek, and Solis computed the determinant of $Th (p,q)$ and found the following formulas.

\begin{introthm}
We have $$ \mathrm{det} (Th (p,q)) = \begin{cases} 
      P_p, & \text{if} \ p \geq 3 \ \text{and} \ q = 2, \\
      L_{2q}-2, & \text{if} \ p = 3 \ \text{and} \ q \geq 3, \\
      S_{p,q}, & \text{if} \ p \geq 4 \ \text{even} \ \text{and} \ q \geq 3, \\
\end{cases}$$ where \begin{itemize}[leftmargin=2em]
    \item $P_p$ is the $p$-th Pell number,\footnote{ \ The \emph{Pell numbers} are defined recursively as $P_1 = 1$, $P_2 = 2$, and $P_k = 2P_{k-1} + P_{k-2}$ for $k \geq 3$.}
    \item $L_{2q}$ is the $(2q)$-th Lucas number,\footnote{ \ Similarly, the \emph{Lucas numbers} are defined recursively as $L_0 = 2, L_1 =1$, and $L_k = L_{k-1} + L_{k-2}$ for $k \geq 2$.}
    \item $S_{p,q}$ is a composite integer given by $$S_{p,q} = q \prod_{1 \leq i \leq q-1, \  1 \leq j \leq \frac{p}{2} -1 } \left ( 4 \sin^2 \left ( \frac{i\pi}{q} \right ) + 4 \sin^2 \left ( \frac{j\pi}{p} \right ) \right ) \cdot $$
\end{itemize}
\end{introthm}

\noindent They also provided a potential composite determinant formula for the remaining case: $p \geq 5$ odd and $q \geq 3$, see \cite[{\sc\S}6]{DM10}.

Mattman and Solis \cite{MS09} later proved the Kauffman--Harary conjecture. The recent work of Bakshi, Guo, Montoya-Vega, Mukherjee, and Przytycki \cite{BGMVMP23} also confirmed the \emph{generalized Kauffman--Harary conjecture}, stating that for every pair of distinct arcs in the reduced alternating diagram of a prime link with determinant $\delta$, there is a Fox $\delta$-coloring that distinguishes these arcs. This conjecture was first proposed by Asaeda, Przytycki and Sikora \cite{APS04}, and proved for Montesinos links.

\subsection{Link Floer and Khovanov Homology Groups} 
\label{sec:floer-khovanov}

We partially describe link Floer and Khovanov homology groups \cite{Kho00, Kho03, OS08_alexander} of the Turk's head link $Th (p,q)$, since such homology theories are well-understood for alternating links by the work of Lee \cite{Lee05}, Ozsv\'ath and Szab\'o \cite{OS08_alexander}, and Ozsv\'ath, Szab\'o and Stipsicz \cite[{\sc\S}10]{OSS15}. 

Given a link $L \subset S^3$ with $\ell$ components, we normalize the Alexander polynomial $\Delta_L (t)$ so that we set $(t^{1/2} - t^{-1/2} )^{\ell - 1} \cdot \Delta_L (t) = \sum_i a_i t^i$. Then we see that $$ \widehat{HFL}_d (Th (p,q) ,s ) = \begin{cases} 
      \Z^{\vert a_s \vert }, & \text{if} \ d = s + \frac{\sigma (Th (p,q) ) - \mathrm{gcd} (p,q) + 1 }{2} \CommaPunct \\
      0, & \text{otherwise.} \\
\end{cases} $$ where $s$ is the Alexander grading. 

Since the link $Th (p,q)$ is fibered, the main result by Ozsv\'ath and Szab\'o \cite{OS05} indicates that we have $$\mathrm{dim} \widehat{HFL}_d ( Th (p,q), s_{top} ) = 1, $$ where $s_{top}$ is the maximal Alexander grading $s \in \Z$ such that the link Floer homology group is non-vanishing. In particular, we have $$s_{top} = g_3 (Th (p,q) ) + \mathrm{gcd} (p,q) - 1 = \frac{(p-1)(q-1) + \mathrm{gcd} (p,q) - 1}{2} \cdot $$ By the work of Ghiggini \cite{Ghi08} and Ni \cite{Ni07}, we further know that link Floer homology detects fibered links in the $3$-sphere, see also Cavallo's article \cite{Cav22}.

For the Turk's head link $Th (p,q)$, the non-vanishing Khovanov homology groups $Kh^{i,j} (Th (p,q))$ lies on the lines: $$j = \begin{cases} 
2i \pm 1, & \text{if} \ p \ \text{is odd}  , \\
2i + q -1 \pm 1, & \text{if} \ p \ \text{is even}. \end
{cases}$$ This is first observed by Mishra and Staffeld as an immediate consequence of Lee's work \cite{Lee05}. 

Therefore, Turk's head links provide examples of \emph{homologically $\sigma$-thin links} in Floer and Khovanov theories, see the article by Manolescu and Ozsv\'ath \cite{MO08} for a further discussion.

\subsection{Unknotting Numbers} 

The \emph{unknotting number} of a knot $K$, denoted by $u(K)$, is defined as the minimum number of crossings that must be changed to unknot it. 

The recent work of Gill, Joshi and Shimizu \cite{GJS24} provided some upper bounds for unknotting numbers of the following specific Turk's head knots. Let $p$ be an odd integer, $n$ be a non-negative integer and $r$ be an integer with $1 \leq r \leq p-1$ and $\gcd(p,r)=1$. Then we have 
$$u(Th (p,np+r)) \leq \frac{n(p^2-1)}{4}+\frac{(p-1)r}{2}-1 . $$ 

\noindent Restricting their attention to more specific cases, they also find \begin{itemize}
    \item for an odd integer $p$ and any positive integer $n$, $$u(Th (p,np+1)) \leq \frac{n(p^2-1)}{4} \CommaPunct$$
    \item for an odd integer $p$ and any non-negative integer $n$, $$u(Th (p,np+2)) \leq \frac{n(p^2-1)}{4} + \frac{p-1}{2} \CommaPunct$$
    \item for an even integer $p$ with $p \geq 4$ and any positive integer $n$ $$u(Th (p,np+1)) \leq \frac{np(p-1)}{2} \cdot$$
\end{itemize}

\begin{figure}[htbp]
\centering
\includegraphics[width=0.3\columnwidth]{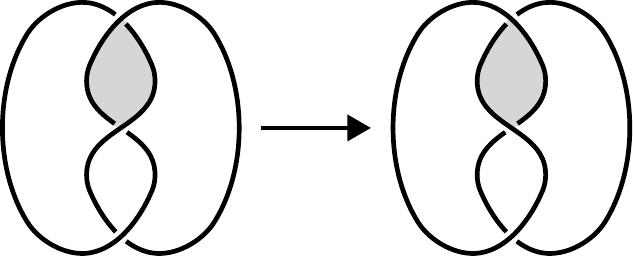}
\caption{A region crossing change.}
\label{fig:region}
\end{figure}

According to Shimizu \cite{Shi14}, Kishimoto introduced the notion of a \emph{region crossing change} on a region $R$ on a knot diagram $D$, which is defined as a local transformation on $D$ obtained by changing the crossings on $\partial R$, see Figure~\ref{fig:region}. Later, Shimizu \cite{Shi14} showed that a region crossing change on a knot diagram is an unknotting operation. Therefore, the \emph{region unknotting number} of a knot diagram $D$, $u_R (D)$, is the minimum number of region crossing changes required to transform $D$ into a diagram of the unknot. In a similar vein, the \emph{region unknotting number} of a knot $K$, $u_R (K)$, is the minimum value of $u_R (D)$ for all minimal crossing diagrams $D$ of $K$.

As observed in \cite[Corollary~5.2]{GJS24}, we, in general, have $$u_R (Th (p,q)) \leq \frac{(p-1)q}{2} + \frac{1}{2} \cdot$$ For the following types of Turk's head knots, Gill, Joshi and Shimizu found sharp upper bounds: \begin{itemize}
    \item for an odd integer $p$ and any positive integer $n$ $$u_R(Th(p,np+1)) \leq \frac{n(p^2-1)}{4} \CommaPunct$$
    \item for an odd integer $p$ and any non-negative integer $n$ $$u_R(Th(p,np+2)) \leq \frac{n(p^2-1)}{4}+\frac{p-1}{2} \cdot$$
\end{itemize}

\begin{figure}[htbp]
\centering
\includegraphics[width=0.35\columnwidth]{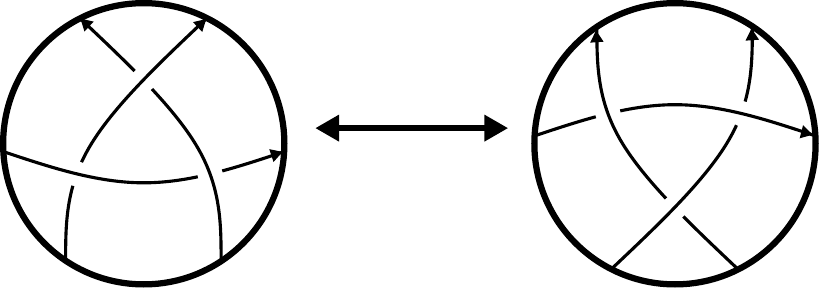}
\caption{The $\Delta$-move.}
\label{fig:delta}
\end{figure}

The $\Delta$-move is a local move for knots, depicted in Figure~\ref{fig:delta}. The \emph{$\Delta$-unknotting number} of a knot $K$ is denoted by $u^\Delta (K) $ and defined as the minimum number of $\Delta$-unknotting operations to change a diagram of $K$ to a diagram for the unknot. See \cite{NY00} for more details.

Nakanishi and Yamada \cite{NY00} computed the $\Delta$-unknotting number of the Turk's head knots $Th (p,q)$. In particular, they found the following formulas:

$$ u^\Delta (Th (p,q) ) = \begin{cases} 
      \sum_{k=1}^{m} Nk^2 = \pm a_2 (Th (p,q) ) , & \text{if} \ p = 2m+1 \ \text{and} \ q = (2m+1)N \pm 1, \\
      \sum_{k=1}^{m} (Nk^2 + k) = - a_2 (Th (p,q) ) , & \text{if} \ p = 2m+1 \ \text{and} \ q = (2m+1)N + 2, \\
      \sum_{k=1}^{m} (Nk^2 - k) =  a_2 (Th (p,q) ) , & \text{if} \ p = 2m+1 \ \text{and} \ q = (2m+1)N -2, \\
      m = \mp a_2 (Th (p,q) ) , & \text{if} \ p = 3m \pm 1 \ \text{and} \ q = 3, \\
\end{cases} $$ where $a_2$ denotes the coefficient of the term $z^2$ in the Conway polynomial.

\subsection{Hyperbolic Structures and Properties}
\label{sec:hyperbolic-general}

Finally, we discuss the rich hyperbolic structures of Turk's head links. A link $L \subset S^3$ is said to be \emph{hyperbolic} if the link complement $S^3 \setminus L$ has a complete Riemannian metric of constant negative curvature and finite volume.

In \cite[p.~758]{MR90}, Meyerhoff and Ruberman highlighted another crucial property of Turk's head links relying on Menasco's work \cite{Men84}, since $Th (p,q)$ is an alternating link. See, for instance, \cite[{\sc\S}11.1]{Pur20} for more details.

\begin{introthm}
For any $p$ and $q$, the Turk's head link $Th (p,q)$ is hyperbolic.
\end{introthm} 

The following definitions relate hyperbolic volumes of knots and links to diagrammatic knot invariants. In \cite{CKP16a}, Champanerkar, Kofman, and Purcell in particular focused on the relationship between the volume $\mathrm{vol}(L)$, the determinant $\mathrm{det}(L)$ and the crossing number $c(L)$ for hyperbolic knots and links $L \subset S^3$. 

A sequence of hyperbolic links $L_n \subset S^3$ is said to be \emph{geometrically maximal} and \emph{diagrammatically maximal} respectively if the following conditions hold: \begin{itemize}
    \item when $c (L_n) \rightarrow \infty$, we have $$ \lim_{n\rightarrow \infty} \frac{ \mathrm{vol} (L_n) }{ c(L_n) } = \upsilon_{\mathrm{oct}} , $$
    \item when $c (L_n) \rightarrow \infty$, we have $$\lim_{n\rightarrow \infty} \frac{ 2 \pi \ \mathrm{log} \ \mathrm{det} (L_n) }{ c(L_n) } = \upsilon_{\mathrm{oct}},$$ 
\end{itemize}

\noindent where $\upsilon_{\mathrm{oct}} \approx 3.66386 $ denotes for the volume of the regular ideal octahedron.

The first definition is related to D. Thurston's well-known observation which provides an upper bound for the hyperbolic volume of $L$ for any diagram of a hyperbolic link $L$: $$\frac{ \mathrm{vol} (L) }{ c(L) } \leq \upsilon_{\mathrm{oct}}.$$ The second definition is inspired by a conjecture of Kenyon for planar graphs \cite{Ken96}, which is known to be equivalent to the following conjecture proposed by Champanerkar, Kofman, and Purcell \cite{CKP16a}: for any link $L$, we have $$\frac{ 2 \pi \ \mathrm{log} \ \mathrm{det} (L) }{ c(L) } \leq \upsilon_{\mathrm{oct}}.$$

Considering the sequence of all the Turk's head links $Th(p,q)$ as $p,q \rightarrow \infty$, Champanerkar, Kofman, and Purcell \cite{CKP16a} proved the following interesting result. 

\begin{introthm}
The sequence of the Turk's head links $Th (p,q)$ is both geometrically and diagrammatically maximal.
\end{introthm}

Let $\upsilon_{\mathrm{tet}} \approx 1.01494 $ denote the volume of the regular ideal tetrahedron. In \cite{CKP16b}, Champanerkar, Kofman, and Purcell provided asymptotically sharp bounds on the volume of Turk's head link: \begin{equation} \label{eq:CKP}
\upsilon_{\mathrm{oct}} (p-2) q \left ( 1- \frac{(2 \pi )^2 }{q^2} \right )^{3/2} \leq  \mathrm{vol} (Th (p,q)) < (\upsilon_{\mathrm{tet}} (p-3) + 4\upsilon_{\mathrm{tet}} )q    
\end{equation} for $p \geq 3$ and $q \geq 7$. They also prove that $S^3 \setminus Th (p,q) $ approaches $\mathbb{R}^3 \setminus \mathcal{W} $ as a geometric limit\footnote{ \ Let $X_n$ and $Y$ be locally compact complete metric spaces. $Y$ is called a \emph{geometric limit} of $X_n$ if there exist basepoints $y \in Y$ and $x_n \in X_n$ such that $(X_n, x_n)$ converges in the pointed bi-Lipschitz topology to $(Y,y)$ as $n \rightarrow \infty$. See \cite{CKP16b} for more details.} as $p,q \rightarrow \infty$, where $\mathcal{W}$ is the infinite alternating weave depicted in Figure \ref{fig:weave} (see \cite{CKP16b}). Note that Burton \cite{Bur16} produced the first examples of a sequence of non-alternating geometrically maximal links by using weaving type tangles.

For a given alternating link diagram, a \emph{visible Conway sphere} intersects the plane of projection in a simple closed curve that meets the link diagram transversely at four points and bounds at least two crossings on each side. In \cite{CKP22}, Champanerkar, Kofman, and Purcell observed that Turk's head links $Th(p,q)$ admit no visible Conway spheres for $p,q \geq3$. 

In \cite{CKP22}, Champanerkar, Kofman, and Purcell associated a set of hyperbolic right-angled ideal polyhedra to any reduced, prime and alternating link diagram. They proved that such a set of right-angled polyhedra is a link invariant, whose volume sum is called the \emph{right-angled volume} of the link $L$ and is denoted by $\mathrm{vol}^\perp (L)$. The initial list of the right-angled volumes of hyperbolic ideal right-angled polyhedra was tabulated by Egorov and Vesnin \cite{EV20}. Later, their tabulation was used to obtain right-angled volume bounds for links by Champanerkar, Kofman, and Purcell \cite{CKP22}, Alexandrov, Bogachev, Egorov and Vesnin \cite{ABEV23}, Egorov and Vesnin \cite{EV24}.

In \cite{CKP22}, Champanerkar, Kofman, and Purcell also proved that $\mathrm{vol}^\perp (L) \leq \mathrm{vol} (L)$ for any hyperbolic alternating link. In particular, they showed that the lower bound is asymptotically sharp by using the Turk's head links. 

\begin{introthm}
We have $$\lim_{p,q \rightarrow \infty} \frac{ \mathrm{vol}^\perp (Th (p,q)) }{ \mathrm{vol} (Th (p,q)) } = 1 .$$
\end{introthm}

An embedded or immersed surface in the knot complement is said to be \emph{totally geodesic} if it is isotopic to a surface that lifts to a set of geodesic planes in the hyperholic space $\mathbb{H}^3$. See \cite{Gan21} for more details. In \cite{Gan21}, Gan proved that all Turk's head links $Th (p,q)$ (except for $Th (3,3)$ and $Th (4,4)$) have both checkerboard surfaces of the standard diagram that are not totally geodesic.

Observe that two checkerboard surfaces of the standard diagram are symmetric, so they become isometric by the Mostow-Prasad rigidity. Therefore, these surfaces are either both totally geodesic or both not totally geodesic.

\begin{figure}[htbp]
\centering
\includegraphics[width=0.45\columnwidth]{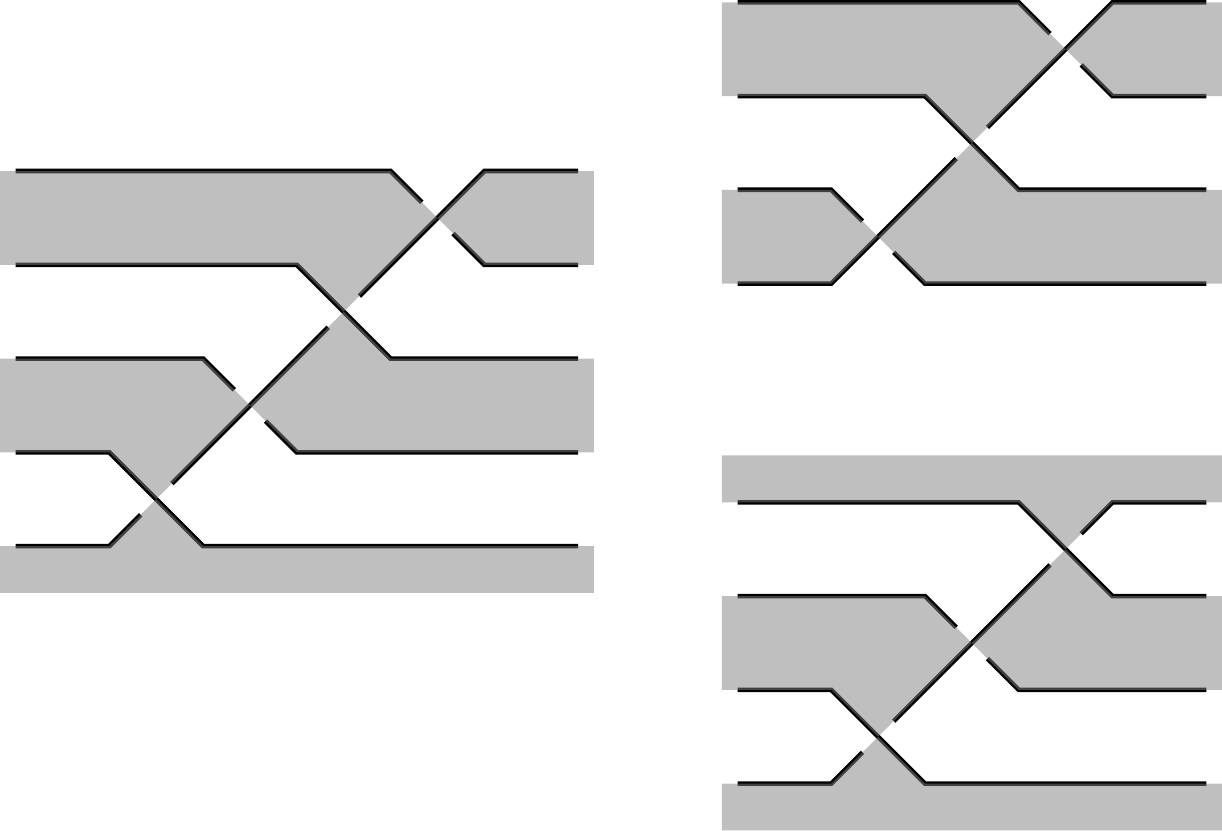}
\caption{The repeating units constituting a checkerboard surface for the Turk's head links and their mirrors. On the left, we have the case for $Th(5,q)$. On the right, we have the case for $Th(4,q)$ (bottom) and $\overline{Th(4,q)}$ (top).}
\label{fig:checkerboard}
\end{figure}

\subsection{Legendrian Properties and Invariants}
\label{sec:legendrian}

A \emph{Legendrian link} $\Lambda$ in a contact $3$-manifold $(M,\xi)$ is a smoothly embedded link which is tangent to the contact planes distribution $\xi$ at every point. It is a well-known fact (see \cite{Etn05, Gei08}) that every smooth isotopy class of links in a contact $3$-manifold $(M,\xi)$ admits a Legendrian representative. In this subsection, we will implicitly refer to Legendrian links as Legendrian links in $S^3$ endowed with the standard contact structure $\xi_{std}$.

Given a Legendrian link $\Lambda$, recall that the \emph{Thurston--Bennequin number} $\tb(\Lambda)$ of $\Lambda$ is defined as the linking number between $\Lambda$ and a nearby pushoff of $\Lambda$ in the positive normal direction of the contact planes. A classical result known as the \emph{Bennequin inequality} \cite{Ben83} states that for a Legendrian link $\Lambda$ the inequality $$\tb(\Lambda)\leq -\chi(\Sigma)$$ holds for any Seifert surface $\Sigma$ for $\Lambda$, where $\chi$ denotes the Euler characteristic.

A well-studied topological invariant for knots and links arising from contact topology is the so-called the \emph{maximal Thurston--Bennequin number}, defined as $$
\TB(L)=\max\{\tb(\Lambda)\;|\;\Lambda\text{ is a Legendrian representative of }L\},
$$where $L\subset S^3$ is a link. Observe from the Bennequin inequality that $\TB(L)$ is always finite.

Since Turk's head links are alternating and non-split, their maximal Thurston--Bennequin numbers are completely determined by either the Jones polynomials and the Trotter--Murasugi signatures \cite[Theorem~4]{Ng05} or by their Kauffman polynomials \cite[Theorem~1.3]{Tan06}.

In particular, by applying \cite[Corollary 1.4]{Tan06} to the diagrams and checkerboard surfaces described in Figure \ref{fig:checkerboard}, we are able to determine their maximal Thurston--Bennequin numbers as follows:
$$
\TB(Th(p,q))=\begin{cases}
    -1-\left(\frac{(p-1)q}{2}\right) \CommaPunct \quad&\quad \text{if} \ p\text{ is odd},\\
    2(q-1)-\frac{pq}{2} \CommaPunct \quad&\quad \text{if} \ p \ \text{is even}.
\end{cases}
$$
As discussed in Section \ref{sec:symmetries}, if $p$ is odd then $Th(p,q)$ is amphichiral, i.e. it coincides with its mirror. In this case, $\TB(Th(p,q)) = \TB(\overline{Th(p,q)})$. However, this is not the case if $p$ is even, and in particular one can compute that 
$$ 
\TB(\overline{Th(p,q)})=-q-\frac{pq}{2} \CommaPunct \quad\quad \text{if} \ p\text{ is even}.
$$ 
Therefore, both $Th(p,q)$ and $\overline{Th(p,q)}$ cannot bound orientable exact Lagrangian surfaces in the standard symplectic $4$-ball $(B^4, \omega_{std})$, (cf.~{\sc\S}\ref{sec:appendix}). This fact is an immediate consequence of Chantraine's main result \cite{Cha10}, see also \cite[Theorem~4.3]{CNS16} in the article by Cornwell, Ng and Sivek. 

\subsection{Sliceness, Surgeries and Fillings}
\label{sec:surgeries}

A link $L \subset S^3$ is said to be \emph{slice} (resp. $\Q$-slice) if $L$ bounds a disjoint union of smoothly properly embedded disks in $B^4$ (resp. in a $\Q$ homology $4$ ball that is a compact connected oriented smooth $4$-manifold $W$ such that $H_* (W; \Q) = H_* (B^4; \Q)$.

For $p$ even, by Theorem \ref{thm:signatures} we know that the signature of $Th(p,q)$ is nonzero, and hence $Th(p,q)$ is not a slice link. Similarly, whenever $Th(p,q)$ is a link, i.e. $\gcd(p,q)>1$, then it is not slice. In fact, if a link $L$ with more than one component is slice, then $\Delta_L\equiv 0$ (see for instance \cite{kawauchi1978alexander}), while by the results in Section \ref{sec:alexander-general} we have that $\Delta_{Th(p,q)}$ is always nonzero.

On the other hand, Conway proved that the knots $Th(3,5)$ and $Th(5,3)$ are slice \cite{Con70}. The sliceness of $Th(7,2)$ was observed by Stoimenow, see KnotInfo \cite{knotinfo}. Note that the recent work of Owens and Swenton \cite{OS24} provided an alternative proof for the sliceness of these three knots.

In {\sc\S}\ref{sec:symmetries}, we also see that the knots $Th (p,q) $ are strongly negative amphichiral if $p$ is odd. Then, by Kawauchi's characterization \cite{Kaw80, Kaw09}, we conclude that $Th (p,q) $ is always $\Q$-slice when $p$ is odd. In this case, we know that $Th(3,q)$ is not slice for certain values of $q$, (see Section \ref{sec:sliceness}) and compare with Problem \ref{prob:sliceness2}.

We call a compact connected oriented smooth $3$-manifold $Y$ a \emph{homology $3$-sphere} if $H_* (Y; \Z) = H_* (S^3; \Z)$. Let $S^3_{p/q} (K)$ be the $3$-manifold obtained by $p/q$ Dehn surgery along $K$ in $S^3$. Relying on Gordon's article \cite{Gor75}, we observe the following $4$-manifold fillings of homology $3$-spheres:  \begin{itemize}
    \item Homology $3$-spheres $S^3_{1/n} (Th(3,5))$, $S^3_{1/n} (Th(5,3))$ and $S^3_{1/n} (Th(7,2))$ bound contractible $4$-manifolds for any $n \in \Z \setminus \{ 0\}$, 
    \item Let $p$ be odd. Then every homology $3$-sphere $S^3_{1/n} (Th(p,q))$ bounds a $\Q$-homology $4$-ball for any $n \in \Z \setminus \{ 0\}$.
\end{itemize}

\noindent For more results in this direction, one can consult \cite[{\sc\S}2.2, {\sc\S}3.1, {\sc\S}3.2]{Sav24} and references therein.

Thurston's hyperbolic Dehn surgery theorem \cite{Thu97} guarantees that for all but finitely many Dehn surgeries along hyperbolic knots yield hyperbolic $3$-manifolds. Therefore, all remaining rational surgeries are said to be \emph{exceptional}. A classification result by Ichihara and Masai \cite{IM16} implies that there is no exceptional surgery for a Turk head knot $Th(p,q)$, except for the figure-eight knot $Th(3,2)$.

\vspace{2em}

\stopcontents[sections]


\section{More Results on the Special Case}
\label{sec:special}

\startcontents[sections]
\setcounter{tocdepth}{2}
\printcontents[sections]{}{1}{}
\vspace{1.5em}

In this section, we consider the subfamily of the $3$-braid Turk's head knots $Th (3,q)$, assuming that $\mathrm{gcd} (3,q) = 1$, which was the subject of a more in-depth analysis in the literature. Using Knotinfo \cite{knotinfo}, Knotscape \cite{knotscape} and SnapPy \cite{snappy}, we identify our knots with a small number of crossings. In particular, we have $$Th (3,2) = 4_1, \ Th (3,4) = 8_{18}, \ Th (3,5) = 10_{123}, \ Th (3,7) = 14_{a19470}, \ \text{and} \ Th (3,8) = 16_{a275159} .$$

\subsection{Fibered Properties}
\label{sec:fibered-properties}

In \cite[Theorem~4.2]{MS21}, Mishra and Staffeld observed that the Turk's head knot $Th (3,q)$ is always fibered by using knot Floer homology $\widehat{HFK}$ \cite{Ghi08, Ni07}. In particular, we have $$ \widehat{HFK}_{d} (Th (3,q), s ) \cong \begin{cases} 
      \Z^{\vert a_s \vert }, & \text{if} \ 0 \leq d \leq q-1, \\
      0, & \text{otherwise,}
\end{cases}$$ so that $$\widehat{HFK}_d (Th (3,q), q-1 ) \cong \Z. $$ This implies that $Th (3,q)$ is fibered (cf.~{\sc\S}\ref{sec:genus-fiberedness}).

An early result of Edmonds and Livingston \cite{EL83} provided a detection result in the sense that the only $q$-periodic fibered knots in $S^3$ of genus $q-1$ for $q \geq 2$ are torus knots $T(3,q)$ and Turk's head knots $Th(3,q)$.

Now we discuss several invariants of hyperbolic fibered links, following Bergbauer's dissertation \cite{Ber98}. Given a fibered link $L$ in $S^3$ with a fiber $F$, it is well-known that its exterior $X_L$ can be described as $$ X_L \doteq S^3 \setminus L = \frac{F \times [0,1]}{(x,1)\sim_h (h(x),0) } $$ for an orientation-preserving diffeomophism $h: F \to F$ with $h \vert_{\partial X_L } = id_{\partial X_L}$. The map $h$ is the so-called \emph{monodromy} of $L$. Moreover, $h$ is said to be \emph{pseudo-Anasov} if there exists a pair of transverse measured singular foliations $(\mathcal{F}^+, \mu^+ )$ and $(\mathcal{F}^-, \mu^- )$ and a real number $\lambda > 1$ such that $$ h (\mathcal{F}^+, \mu^+ ) = (\mathcal{F}^+, \lambda \mu^+ ) \ \ \text{and} \ \ h (\mathcal{F}^-, \mu^- ) = (\mathcal{F}^-, \lambda^{-1} \mu^- ) .$$

\noindent The number $\lambda$ is called the \emph{stretching factor} or \emph{dilatation}. Due to Thurston \cite{Thu97}, the monodromy of a hyperbolic fibered link is freely isotopic to a pseudo-Anosov homeomorphism $\varphi$. Note that $\varphi \vert_{\partial X_L } $ is never equal to $id_{\partial X_L}$. Moreover, the stable lamination on $F$ suspends to a lamination $\mathcal{L}$ on $X_L$ of codimension $1$. The cusps of the complementary region of $\mathcal{L}$ containing $\partial X_L$ are parallel to $n$ copies of an essential simple closed curve $\alpha = (a,b) \subset \partial X_L$. In this case, $n \alpha$ is called the \emph{degeneracy slope} $d(L)$ of the link $L$, see also \cite[Definition~5]{GO89}.

Bergbauer computed the three invariants of fibered links for the Turk's head link $Th(3,q)$ as follows.

\begin{introthm}
We have that \begin{itemize}
    \item the monodromy of $Th (3,q)$ is $h = R^{q-1} \circ L^{q-1} $, see \cite[p.~39]{Ber98}. Here, $R$ and $L$ denote the monodromies of the positive Hopf link $T(2,2)$ and the negative Hopf link $\overline{T(2,2)}$, respectively,
    \item the stretching factor $Th (3,q)$ is $\lambda = \frac{3 + \sqrt{5} }{2}$ that is the square of the golden ratio, see \cite[Theorem~3.2.1]{Ber98},
    \item the degeneracy slope of $Th (3,q)$ is $d(Th (3,q)) = q (1,0)$, see \cite[Theorem~3.2.2]{Ber98}.
\end{itemize}
\end{introthm}

With the result of Gabai and Oertel \cite{GO89}, the computation of these invariants has an immediate consequence. Bergbauer showed that the universal cover of a $3$-manifold obtained by any non-trivial Dehn surgery on a Turk's head knot is $\mathbb{R}^3$.

In \cite{Lie17}, Liechti studied the minimal stretching factors arising from Penner's construction \cite{Pen88} on closed orientable surfaces, and observed that the minimum among all stretching factors arising from Penner’s construction is geometrically realized by the monodromy of the figure eight knot $Th(3,2)$. By Bergbauer's result, the minimum is also realized by the other Turk's head knots $Th(3,q)$.

\subsection{Hyperbolicity and Symmetry Groups}
\label{sec:symgroups}

In \cite[Example~6.8.11]{Thu97}, Thurston observed that the Turk's head knot $Th (3,q)$ is hyperbolic for any $q$. This was also reproved by Long \cite[{\sc\S}7]{Lon83}. One can check the paper by Thistlethwaite and Tsvietkova \cite{TT14} for a new alternative proof.

In \cite[{\sc\S}5.2]{CKP22}, the right-angled volumes $\mathrm{vol}^\perp (Th (3,q))$ were explicitly computed for $q \leq 7$. Champanerkar, Kofman, and Purcell also observed that the computation of $\mathrm{vol}^\perp (Th (3,q))$ is equivalent to solving a one-variable polynomial equation.

The \emph{symmetry group} of a knot $K$ in $S^3$, which is denoted by $\mathrm{Sym} (S^3, K)$, is defined to be the mapping class group of the knot exterior $S^3 \setminus \nu(K)$, see Kawauchi's book \cite[{\sc\S}10.6]{Kaw90}. By the work of Sakuma and Weeks \cite[Proposition~I.2.5]{SW95}, the symmetry group of the Turk's head knot $Th (3,q)$ is completely determined. 

\begin{introthm}
For any $q$, we have $$\mathrm{Sym} (S^3 , Th (3,q)) \cong D_{2q},$$ where $D_{n}\cong\Z/n\Z\rtimes\Z/2\Z$ denotes the dihedral group of $2n$ elements.
\end{introthm}

\subsection{Alexander and Conway Polynomials}
\label{sec:alexander-conway}

Using skein relations for Conway polynomials, Takemura computed in \cite[Theorem~4.3]{Tak16} that we have the following recursive formula $$\nabla_{Th (p,3)} (z) = \sum_{i = 0}^\infty a_i z^i \ \ \text{and} \ \ \nabla_{Th (3,p)} (z) = \sum_{i = 0}^\infty b_i z^i $$ where $$\begin{cases} 
      a_i = b_i , & \text{for} \ i=0,1 \ (\text{mod} \ 4),\\
      a_i = -b_i , & \text{for} \ i=2,3 \ (\text{mod} \ 4).  
\end{cases}$$ One can find more explicit formulations for the Conway polynomials by relating them with Alexander polynomials. 

Recently, AlSukaiti and Chbili calculated Alexander polynomials of Turk's head knots $Th(3,q)$ explicitly \cite[Theorem~1.2]{AC23}, and found the following closed formula. In \cite[Proposition~5.1]{AC23}, they also computed the roots of the Alexander polynomial $\Delta_{Th(3,q)} (t)$.

\begin{introthm}
We have $$\Delta_{Th (3,q)} (t) = \sum_{k=0}^{2q-2} \alpha_{q,k} t^{k}, \ \ \alpha_{q,k} = \sum_{i=0}^{q-1} (-1)^{k-q+i+1} \frac{2q(q+i)!}{(q-i-1)! (2i+2)!}  \begin{pmatrix}
			i;3 \\
			k+i-q+1
		\end{pmatrix}, $$ where $\begin{pmatrix}
			q;3 \\ k
		\end{pmatrix}$ is the coefficient of $z^k$
in the expansion of $(1+z+z^2)^q$. Moreover, the roots of $\Delta_{Th(3,q)} (t)$ are of the form: $$z = -\frac{1}{2} \left( 2\cos \left( \frac{2k}{n} \pi \right ) -1 \pm \sqrt{ \left (2\cos \left (\frac{2k}{n} \pi \right ) -1 \right ) ^2 -4 } \right) \CommaPunct $$ with $1 \leq k \leq \lfloor n/2 \rfloor.$
\end{introthm}

As an application, they showed that the Turk's head knots $Th (3,q)$ satisfy the \emph{Fox’s trapezoidal conjecture} \cite{Fox61}: given an alternating knot $K$ with the Alexander polynomial $\Delta_K(t)=\pm \displaystyle\sum_{i=0}^{2n}\alpha_i(-t)^i$, with $\alpha_i>0$, there exists an integer $r \leq n$ such that $$\alpha_0<\alpha_1< \dots < \alpha_{n-r}=\dots =\alpha_{n+r}> \dots ...>\alpha_{2n-1}>\alpha_{2n}.$$ They also observed that the zeroes of $\Delta_{Th (3,q)} (t)$ satisfy the \emph{Hoste's condition}: $$\mathfrak{Re}(z) > -1 .$$ Moreover, if $z$ is a non-real root then $|z| =1$. Note that certain alternating knots, which provide counter-examples to Hoste's condition were found by Hirasawa, Ishikawa and Suzuki \cite{HIS19}.

\subsection{Jones Polynomials}
\label{sec:jones}

AlSukaiti and Chbili also calculated the Jones polynomial of Turk's head knots $Th (3,q)$ explicitly \cite[Corollary~3.4]{AC23}, obtaining the following formulation. 

\begin{introthm}
We have $$ V_{Th (3,q)}(s) =\sum_{k=0}^{2q} a_k s^{k-q},$$ where $$a_k = \begin{cases}
c_{q,k}, & \text{if} \ |k-q| \neq 1,\\
c_{q,k}-1 , & \text{if} \ |k-q| = 1,
\end{cases} \quad  \text{and} \quad c_{q,k} = \begin{cases}
		\sum_{j=0}^{\lfloor k/2 \rfloor} \frac{q}{q-j}\binom{q-j}{q-k+j}\binom{k-j-1}{j} \CommaPunct & \text{if} \ k < 2q, \\
		1, & \text{if} \ k = 2q.
  \end{cases}$$
\end{introthm}
  
\noindent Note that various sample or full computations for Jones polynomials of $Th (3,q)$ previously appeared in the articles by Mishra and Staffeld \cite{MS21} and Joshi, Negi and Prabhakar \cite{JNP23}. For the sample computations of HOMFLY-PT polynomials, colored HOMFLY-PT polynomials and colored Jones polynomials of $Th (3,q)$, one can consult the paper of Singh and Chbili \cite{SC22}, Singh, Mishra, and Ramadevi \cite{SMR21} and Kaiser and Mishra \cite{KM24}, respectively.

\subsection{Determinants} 
\label{sec:determinants2}

Using the Alexander or Jones polynomial of computation of Mishra and Staffeld, we can obtain the determinant formula for $Th (3,q)$ \cite[Corollary~3.5]{AC23} as follows 
\begin{equation}
\label{eq:determinant}
\mathrm{det}(Th (3,q)) = \vert \Delta_{Th (3,q)} (-1) \vert = \vert V_{Th (3,q)} (-1) \vert =  L_{2q} - 2 = F_{2q-1} + F_{2q+1} - 2 = L_q F_q -2
\end{equation}
\noindent where $F_k$ and $L_k$ denote the $k^{th}$ Fibonacci and Lucas numbers\footnote{ \ The \emph{Fibonacci numbers} are defined by the recurrence relation by $F_0 = 0, F_1 =1$, and $F_k = F_{k-1} + F_{k-2}$ for $k \geq 2$. Recall that the \emph{Lucas numbers} are defined recursively as $L_0 = 2, L_1 =1$, and $L_k = L_{k-1} + L_{k-2}$ for $k \geq 2$.}, respectively. 

We want to emphasize that $\mathrm{det}(Th (3,q))$ was previously computed in several knot theory articles with different methods. See the papers by Stoimenow \cite[Lemma~6.13]{Sto05} and \cite[Lemma~3.1]{Sto07}, Przytycki \cite[Example~3]{Prz11}, Dowdall, Mattman, Meek, and Solis \cite[Theorem~4]{DM10}, Lopes and Matias \cite[Corollary~1.1]{LM15}, and Joshi, Negi and Prabhakar \cite[Theorem~2.5]{JNP23}. 

Since the number of spanning tress of a wheel graph is equal to the determinant of $Th(3,q)$, the formula in (\ref{eq:determinant}) was found in various graph theory articles in the late 1960s; see, for instance, the papers by Sedl{\'a}{\v c}ek \cite{Sed69}, Myers \cite{Mye71, Mye75}, Harary, O'Neil, Read and Schwenk \cite{HORS72}, and Hilton \cite{Hil72}.

\subsection{Colorings} 

Let $\# col_n Th(3, q)$ denote the number of $n$-colorings of $Th(3, q)$. In \cite{LM15}, Lopes and Matias obtained the following formula for the colorings: $$ \# col_n Th(3, q) =
\begin{cases}
n \cdot ( \mathrm{gcd} (u_{q-1}, n) )^2 , & \text{ if $q$ is odd,}\\
n \cdot \mathrm{gcd} (5u_{q-1}, n) \cdot \mathrm{gcd} (u\sb{q-1}, n), & \text{ if $q$ is even,}
\end{cases} $$ where $$u\sb{q}=\frac{1}{\sqrt{5}}\Bigg[ \Bigg( \frac{1+\sqrt{5}}{2} \Bigg)\sp{q+2} - \Bigg( \frac{-1 + \sqrt{5}}{2} \Bigg)\sp{q}
 - \Bigg( \frac{1-\sqrt{5}}{2} \Bigg)\sp{q+2} + \Bigg( \frac{-1-\sqrt{5}}{2} \Bigg)\sp{q}\Bigg] .$$

\vspace{0.5em}
 
\noindent Restricting on the specific cases, they found more explicit formulas as follows. Given $n, r \in \Z^{>0}$, we have

\begin{enumerate}
\item $2\mid n \text{  and } 3\mid q  \qquad \text{ if and only if }\qquad   mincol\sb{n}Th (3, q) = 2$,

\item $\Big( 2\nmid n  \text{  or }   3\nmid q\Big)\,   \text{  and }  \, 3\mid n \, \text{  and } \, 4\mid q \quad \text{ if and only if } \quad mincol\sb{n}Th(3, q) = 3$,
  
\item
\[
\Big( 2\nmid n  \text{ or }  3\nmid q \Big) \, \text{ and } \, \Big( 3\nmid n  \text{ or }  4\nmid q \Big)  \, \text{ and }  \,  \Big[\big( \, 5\mid n \text{ and } 2\mid q\, \big) \text{ or } \big(\,  7\mid n \, \text{ and }\,  8\mid q\, \big) \Big] \qquad
\]
\[
\text{ if and only if  }\qquad mincol\sb{n}Th(3, q) = 4,
\]
\item If
\[
\Big( 2\nmid n \, \text{ or }  \, 3\nmid q \Big)  \, \text{ and }  \, \Big( 3\nmid n   \, \text{ or } \, 4\nmid q \Big)  \, \text{ and }  \, \Big( 5\nmid n  \, \text{  or }  \, 2\nmid q\Big)  \, \text {and }  \, \Big( 7\nmid n  \, \text{ or }  \, 8\nmid q\Big)
\]
\[
\text{and  } \qquad \Big( 11\mid n   \text{ and } 5\mid q \Big)  \qquad \text{ then } \qquad mincol\sb{n}Th(3, q) = 5.
\]
\end{enumerate} Here, $$mincol\sb{n} K = min \{ n_{D_K , n} \ \vert \ D_K \ \text{is a diagram of a knot} \ K \}$$ where $ n_{D_K , n}$ is the least number of colors it takes to set up a non-trivial $n$-coloring of $D_K$.

In \cite{CGP23}, Christiana, Guo and Przytycki studied the \emph{reduced group of Fox colorings} $Col^{red}(D)$ for Turk's head knots $Th (3,q)$. Recall that $Col^{red}(D)$ is defined as $$\displaystyle Col^{red}(D) \doteq \frac{Col(D)}{Col^{trivial}(D)}$$ where $Col(D)$ is the abelian group whose generators are indexed by the arcs of a diagram $D$ of a knot, and $Col^{trivial}(D)$ the group of trivial colorings of $D$. Remark that $Col^{red}(D)$ can be interpreted as the first homology of the double branched covering of $S^3$ branched along $D$. See \cite[{\sc\S}1]{CGP23} for more details.

Using the standard reduced alternating diagram $D_q$ (see Figure~\ref{fig:ourknots}) for $Th (3,q)$ (which is also known as the wheel graph, see \cite[Figure~3]{CGP23}), Christiana, Guo and Przytycki computed their reduced groups of Fox colorings completely: $$Col^{red}(D_q)=
\left\{
\begin{array}{ll}
\mathbb Z_{F_{q-1}+F_{q+1}} \oplus \mathbb Z_{F_{q-1}+F_{q+1}}, & \mbox{if $q$ is odd},\\
\mathbb Z_{5F_{q}} \oplus \mathbb Z_{F_{q}}, & \mbox{if $q$ is even,}
\end{array} 
\right.$$ where $F_k$ is the $k^{th}$ Fibonacci number.

\subsection{Hyperbolic Structures and Properties}
\label{sec:hyperbolic}

In this subsection, we no longer assume that $\mathrm{gcd} (3,q) =1 $, so we rather consider the Turk's head links $Th (3,q)$.

A link $L$ is called ($2\pi /n$)\emph{-hyperbolic} if its $n$-fold cyclic branched covering is a hyperbolic manifold and the covering group acts by isometries. In \cite[{\sc\S}4.3]{BS79}, Bonahon and Siebenmann showed that $Th (3,q)$ for $q \geq 4$ is $\pi$-hyperbolic. In another direction, Paoluzzi \cite{Pao99} proved that if a knot is fibered and $\pi$-hyperbolic then it is determined by its $2$-fold cyclic branched covering together with its $n$-fold cyclic branched coverings for any $n > 2$. 

Given a $3$-manifold $M$, the \emph{virtual first Betti number} of $M$ is the supremum of the first Betti number of all finite coverings. In \cite{KL88}, Kojima and Long showed that the virtual first Betti number of the manifold obtained by $\mu / \lambda$ Dehn surgery on the Turk's head knot $Th (3,2q)$ is infinite, provided $\mu = 0 \mod{4\lambda}$ and $\mu/\lambda \neq \pm 8q$. This can be seen as early support to Waldhausen's \emph{virtual Haken conjecture}, which was later proved by Agol \cite{Ago13}. For a detailed discussion, see the book by Aschenbrenner, Friedl and Wilton \cite[{\sc\S}5.9]{AFW15}.

Meyerhoff and Ruberman \cite{MR90} also studied the Atiyah-Patodi-Singer $\eta$-invariant by using the Turk's head knots $Th (3,q)$. Let $M$ be the hyperbolic $3$-manifold obtained by $0$-surgery on $Th (3,q)$, and let $M^\varphi$ be the $3$-manifold obtained by cutting and pasting by $q$-periodic symmetry $\varphi$: $$ M^\varphi \doteq \left ( M \setminus \left ( F \times \left \{ \frac{1}{2} \right \} \right) \right) \cup \left ( F \times (0,1) \right), $$ where the gluing is via the identity on $F \times (0,\frac{1}{2})$, and via $\varphi$ on $F \times (\frac{1}{2}, 1)$. Then, they obtained the following formula for the the Atiyah-Patodi-Singer $\eta$-invariant of $M^\varphi$: $$ \eta (M^\varphi) = \begin{cases} 
      \frac{-8 (q-1)^2 }{9q} \pm 1, & \text{if} \ q = 1  \ (\text{mod} \ 3)  , \\
      \frac{ 4 (q-2) (1-2q) }{9q} \pm -1, & \text{if} \ q = -1  \ (\text{mod} \ 3) . 
\end{cases}$$ Using these manifolds, they also showed that any rational number is realized as the Chern-Simons invariant of a compact hyperbolic $3$-manifold.

Now, we present the results regarding Fibonacci manifolds and Turk's head links $Th (3,q)$. Following Helling, Kim and Mennicke \cite{HKM98}, the \emph{Fibonacci manifolds} $M_q$ for $q \geq 2$ are defined as the following quotient $$ M_q = X_q / F(2,2q)$$ where $X_2$ is the $3$-sphere $S^3$, $X_3$ is the Euclidean space $\mathbb{R}^3$, $X_q$ is the hyperbolic space $\mathbb{H}^3$ for $q \geq 4$, and $F(2,2q)$ is the \emph{Fibonacci group} introduced by Conway in the mid-1960s: $$ F(2,2q) \doteq \langle x_1, x_2, \ldots, x_q \ \vert \ x_i x_{i+1} = x_{i+2}, \quad i \ \text{mod} \ q \rangle .$$

Hilden, Lozano and Montesinos \cite{HLM92} proved that every Fibonacci manifold $M_q$ can be represented as a $q$-fold branched covering of $S^3$ branched over the figure-eight knot $Th (3,2)$. We can also represent $M_q$ as the $2$-fold branched covering of $S^3$ branched over the link $Th (3,q)$, see the article by Mednykh and Vesnin \cite{MV96}. We can also relate the volumes of $M_{2q}$ with the volumes of the link complements $S^3 \setminus Th (3,q)$ in the following fashion: $$\mathrm{vol} (M_{2q}) = \mathrm{vol} (S^3 \setminus Th (3,q)) = 4q ( \Lambda(\alpha + \gamma) + \Lambda(\alpha - \gamma) ), $$ where $$ \Lambda (x) = - \int_{0}^x \ln \vert 2 \sin \zeta \vert d \zeta $$ is the Lobachevsky function, $\gamma = \frac{\pi}{2q}$ and $\alpha = \frac{1}{2} \arccos \left (\cos (2\gamma) - \frac{1}{2} \right )$. In particular, we have the asymptotic formula for the volume as $q \rightarrow \infty$: $$ \mathrm{vol} (S^3 \setminus Th (3,q)) \rightarrow 8q \ \Lambda\left(\frac{\pi}{6} \right ) \cdot$$ These results were also proved by Mednykh and Vesnin \cite{MV95a, MV95b} which relied on the work of Thurston \cite{Thu97}. For more details, one can see the comprehensive survey by Cavicchioli, Repov{\v s} and Vesnin \cite{CRV01}. A follow-up result by Matveev, Petronio and Vesnin \cite{MPV09} also indicated that for large $q$ the following inequalities hold:
\begin{equation}
\label{ineq:fibonacci}
2q - \frac{34}{q} < \frac{\mathrm{vol} (S^3 \setminus Th (3,q))}{v_3} < 2q.
\end{equation}

In \cite{DL07}, Dasbach and Lin introduced higher twist numbers of a knot in terms of the Jones polynomial which are related to the volume of a knot. Given a knot $K$, let $$ V_K (s) = a_{-m} s^{-m} + a_{-m+1} s^{-m+1} + \ldots + a_{n-1} s^{n} + a_{n} s^{n} $$ so that $j$th \emph{twist number} of $K$ is defined as $$\mathrm{Tw}_j (K) = \vert a_{-m + j} \vert + \vert a_{-n - j} \vert .$$

Mishra and Staffeldt \cite{MS21} computed the second and the third twist number of the Turk's head knots $Th (3,q) $ as follows: \begin{itemize}
    \item $\mathrm{Tw}_2 (Th (3,q)) = q (q-1) $ for $q \geq 5$,
    \item $\mathrm{Tw}_3 (Th (3,q)) = \frac{q (q-1) (q-2) }{3} + 2q$ for $q \geq 5$.
\end{itemize}

\noindent They compared these computations with the volumes of the Turk's head knots $Th (3,q)$ and found better lower bounds compared to the lower bounds in (\ref{eq:CKP}). See \cite[{\sc\S}7]{MS21} for further details. With the Jones polynomial computations in {\sc\S}\ref{sec:jones}, such relations might be well-understood in the near future. Focusing on certain examples of Turk's head knots in \cite[{\sc\S}6]{MS21}, Mishra and Staffeldt observed an interesting phenomenon: the ranks of the Khovanov homology groups compared with certain their density functions roughly approximate a normal distribution. See also the paper by Mishra and Raundal \cite{MR19} for a comprehensive discussion.

Finally, we briefly mention the work of Hilden, Lozano and Montesinos \cite{HLM00}, in which they studied character varieties of Turk's head links $Th (3,q)$ for $q \leq 5$. They computed the excellent component of the character curve, which is the component that contains the complete hyperbolic structure and can be used to study volume, Chern-Simons invariant and arithmetic properties.

\subsection{Sliceness and Concordance}
\label{sec:sliceness}

As an immediate consequence of the Fox-Milnor condition \cite{FM66}, the determinant of a knot must be a square integer. Using the identity (\ref{eq:determinant}), we can see that the Turk's head knots $Th (3,q)$ are not slice for $q$ even: $$\det(Th (3,q) )=L_{2q}-2=(L_q)^2-2\cdot(1+(-1)^{q}).$$ 

When $q$ is odd, the Turk's head knots $Th (3,q)$ are both strongly positive amphichiral and strongly negative amphichiral. In this case, various classical and new obstructions for sliceness automatically vanish. In particular, the following invariants do not provide sliceness obstructions: $\tau$-invariant \cite{OS03}, $\delta$-invariant \cite{MO07}, $V_0$-invariant \cite{OS08_surgeries}, $s$-invariant \cite{Ras10}, $\nu$-invariant \cite{OS11}, $s_N$-invariants \cite{Lob09, Lew14}, $\delta_p$-invariants \cite{Jab12}, $s^\sharp$-invariant \cite{KM13}, $\epsilon$-invariant \cite{Hom14}, $\nu^+$-invariant \cite{HW16}, $\Upsilon$-invariant \cite{OSS17}, $\underline{V}_0$- and $\overline{V}_0$-invariants \cite{HM17}, $\Upsilon^2$-invariant \cite{KL18}, $\mathrm{HFB}_\mathrm{conn}$-invariant \cite{AKS20}, $\tau^\sharp$-invariant \cite{BS21}, $\nu^\sharp$-invariant \cite{BS21}, $\varphi_j$-invariants \cite{DHST21}, $\theta^q$-invariants \cite{Bar22}, $\tilde{s}$-invariant \cite{DISST22}, and $\tilde{\epsilon}$-invariant \cite{DISST22}.

On the one hand, it is well-known that the Turk's head knot $Th (3,5)$ is slice, see for instance \cite{Con70}. On the other hand, in his dissertation \cite{Lon83}, Long proved that the Turk's head knot $Th (3,7)$ is not \emph{homotopy-ribbon}\footnote{A homotopy-ribbon knot is a knot bounding a slice disk such that the inclusion map from the knot exterior to the slice disk exterior is $\pi_1$-surjective.}.

Despite the difficulties, however, there are positive results for sliceness obstructions based on the powerful theory of twisted Alexander polynomials\footnote{The twisted Alexander polynomials were introduced in the influential work of X. S. Lin \cite{Lin01}. For further details, see the survey article by Friedl and Vidussi \cite{FV11}.} due to Herald, Kirk and Livingston \cite{HKL10}. Sartori improved Long's result in his master's thesis \cite{Sar10} by showing that $Th (3,7)$ is not slice. It was recently highlighted by Lamm \cite{Lam23} that $Th (3,7)$ is the minimum crossing knot that is prime strongly positive amphichiral but not slice. Furthermore, Aceto, Meier, A. Miller, M. Miller, Park and Stipsicz proved that the Turk's head knots $Th (3,11)$, $Th (3,17)$ and $Th (3,23)$ are not slice, see \cite{AMMMPS21}. 

\begin{introthm}
The Turk's head knots $Th (3,7)$, $Th (3,11)$, $Th (3,17)$ and $Th (3,23)$ are not slice.
\end{introthm}

We can expand this discussion by introducing the knot concordance groups. Two knots $K_0$ and $K_1$ in $S^3$ are called \emph{concordant} (resp. $\Q$\emph{-concordant}) if there exists a properly and smoothly embedded annulus $A$ in $S^3 \times [0,1]$ (resp. in a $\Q$-homology $S^3 \times [0,1]$) such that $\partial A = -(K_0) \sqcup K_1$, i.e., the knots $K_0$ and $K_1$ cobound the annulus $A$. Both concordance and $\Q$-concordance yield equivalence relations on the set of oriented knots so that $$ \mathcal{C} = \frac{\{ \text{oriented knots in} \ S^3 \}}{\text{concordance}} \ \ \text{and} \ \ \mathcal{C}_\Q = \frac{\{ \text{oriented knots in} \ S^3 \}}{\Q\text{-concordance}}$$ admit abelian group structures with operations induced by connected sum. The groups $\mathcal{C}$ and $\mathcal{C}_\Q$ are respectively called the \emph{concordance group} and \emph{$\Q$-concordance group}. A classical observation indicates that a knot is slice (resp. $\Q$-slice) if and only if it is concordant to the unknot in $\mathcal{C}$ (resp. in $\mathcal{C}_\Q$). See the survey articles \cite{Liv05, Hom17, Sav24} for further discussions.

Recall that all the Turk's head knots $Th (3,q)$ are strongly negative amphichiral. This has two important implications: \begin{enumerate}
    \item The order of a Turk's head knot $Th (3,q)$ is at most $2$ in $\mathcal{C}$,
    \item The order of any Turk's head knot $Th (3,q)$ is $1$ in $\mathcal{C}_\Q$, since any strongly negative amphichiral knot is $\Q$-slice due to Kawauchi's characterization \cite{Kaw80, Kaw09}.
\end{enumerate}

\noindent Therefore, we can conclude from our above discussion that \begin{itemize}
    \item The order of $Th (3,5)$ is $1$ in $\C$,
    \item The order of $Th (3,q)$ is $2$ if $q$ is even or $q= 7, 11, 17, 23$,
    \item For the remaining infinitely many cases, the exact order of $Th (3,q)$ is unknown.
\end{itemize}

We can also completely describe their orders in the \emph{algebraic concordance group} $\mathcal{AC}$, introduced by Levine \cite{Lev69a, Lev69b}. Recall that two Seifert matrices $V_1$ and $V_2$ are called \emph{algebraically concordant} if $V_1 \oplus -V_2$ is metabolic, i.e., $A \oplus -B$ is congruent to a bloc matrix of the form $\begin{pmatrix}
    0 & B \\
    C & D
\end{pmatrix} \cdot$ The set of Seifert matrices modulo the algebraic concordance relation forms the abelian group called algebraic concordance group $\mathcal{AC}$.\footnote{ \ In \cite{Lev69a, Lev69b}, Levine completely described the structure of $\mathcal{AC}$ by showing that $\mathcal{AC} \cong \Z^\infty \oplus \Z_2^\infty \oplus \Z_4^\infty$. There is also a surjective group homomorphism $\phi : \mathcal{C} \to \mathcal{AC}$, which is not injective by the work of Casson and Gordon \cite{CG78,CG86}. See Livingston's survey \cite[{\sc\S}2.2]{Liv05} for more details.} 

The Turk's head knots $Th (3,q)$ are strongly positive amphichiral for odd values of $q$, so they are algebraically slice and have order $1$ in $\mathcal{AC}$ by the classical result of Long \cite{Lon84}. The remaining Turk's head knots are shown to be non-slice and have order $2$ in $\mathcal{AC}$ by Brandenbursky \cite{Bra16} in which their Arf invariants were computed. This result was independently proved by Lisca \cite{Lis17} where the sliceness obstructions are obtained by using Donaldon's diagonalization theorem \cite{Don83, Don87}. Both proofs also recover the above fact that the order of $Th (3,q)$ is $2$ in $\mathcal{C}$ when $q$ is even.
\stopcontents[sections]


\section{Open Problems}
\label{sec:problems}

Recall that we assume $p \geq 3$ and $q \geq 2$ for Turk's head links $Th (p,q)$, and we call them Turk's head knots if $\mathrm{gcd} (p,q) =1$.

\subsection{On Diagrammatic Properties}

In {\sc\S}\ref{sec:crossing}, we explained that most of the Turk's head knots are not perfectly straight. We pose the following problem to determine the status of being perfectly straight for the rest of the Turk's head knots $Th(p,q)$.

\begin{introprob}
Determine whether the Turk's head knot $Th (p,q)$ is perfectly straight when $p \geq 3$ and $q < p$.
\end{introprob}

Consider a minimal diagram $D$ of an alternating knot $K$ with $c(K) = n$. If changing a crossing in $D$ gives a minimal diagram $D'$ of a knot $K'$ with $c(K') = n$, then $K$ is said to have a \emph{Kauffman property}. An alternating knot $K$ admitting a minimal diagram $D$ such that all crossings of $D$ have the Kauffman property is called a \emph{super-Kauffman knot}.

In \cite{JKSSZ16}, Jablan, Kauffman, Sazdanovic, Sto{\v s}i{\'c}, and Zekovi{\'c} proposed an interesting conjecture on crossing changes in minimal alternating diagrams by focusing on Turk's head knots.

\begin{introconj}
Prove that all Turk's head knots $Th(p,q)$ for $p \geq 5$ and for $q \geq 4$ are super-Kauffman knots.
\end{introconj}

\subsection{On Fibered Properties}

In \cite{GK90}, Gabai and Kazez conjectured that the degeneracy slope of a non-torus fibered alternating knot is of the form $d(K) = n(1,0)$ for some integer $n >1$. We highlight their conjecture by addressing Turk's head knots:

\begin{introprob}
    Verify that $d(Th(p,q)) = n(1,0)$ for some integer $n >1$. If this is the case, how does $n$ depend on $p$ and $q$?
\end{introprob}

In \cite[p.~33]{Ber98}, Bergbauer noted that for a fixed $p$, all the Turk's head knots $Th (p,q)$ have the same stretching factor. See also the discussion on Penner's construction in {\sc\S}\ref{sec:fibered-properties}. For the other fundamental invariants of fibered links, we raise the following problem.

\begin{introprob}
    Compute the monodromies and stretching factors of all Turk's head links $Th (p,q)$.
\end{introprob}

\subsection{On Symmetries and Symmetry Groups}

In {\sc\S}\ref{sec:symmetries}, we saw that $Th (p,q) $ is always periodically positive amphichiral if $p$ is odd. The following problem is about their characterization in terms of the order of given maps. For the first part of the problem, Lamm's SnapPy script \cite[{\sc\S}2]{Lam23} would be helpful.

\begin{introprob}
Let $p$ be odd. Determine which Turk's head links $Th (p,q)$ are strongly positive amphichiral. Moreover, determine the orders of the maps realizing periodic positive amphichirality.
\end{introprob}

As we discussed in {\sc\S}\ref{sec:symgroups}, we know that $$\mathrm{Sym} (S^3 , Th (3,q)) \cong D_{2q} \ \ \text{for a fixed value of} \ q.$$ We expect the following conjecture to be true, relying on the experiments performed by SnapPy \cite{snappy}. When $p$ and $q$ are both odd, a result of Kodama and Sakuma \cite[{\sc\S}1]{KS92} implies that its symmetry group is isomorphic to $D_{2n}$ for some $n$.

\begin{introconj} Suppose $\gcd(p,q)=1$. Then the Turk's head knot $Th (p,q)$, has symmetry group given by $$\mathrm{Sym} (S^3 , Th (p,q)) \cong D_{2q}.$$
\end{introconj}

\subsection{On Polynomial and Homological Invariants}

In {\sc\S}\ref{sec:alexander-general}, {\sc\S}\ref{sec:determinants}, {\sc\S}\ref{sec:floer-khovanov}, {\sc\S}\ref{sec:alexander-conway}, and {\sc\S}\ref{sec:jones}, we shared various results on the computations of polynomial and homological invariants of Turk's head links. To understand their properties better, we pose the following general problems.

\begin{introprob}
\label{prob:polynomial}
Compute the following invariants of $Th (p,q)$ in terms of $p$ and $q$ more explicitly:
\begin{itemize}
    \item Alexander polynomial,
    \item Conway polynomial,
    \item Jones polynomial.
\end{itemize}
\end{introprob}

\vspace{0.5 em}

\noindent Note that Takemura conjectured in \cite[{\sc\S}5]{Tak16} that the Conway polynomial $\nabla_{Th (p,q)} (z)$ is of the following form: $$\nabla_{Th (p,q)} (z) = \sum_{i = 0}^\infty a_i z^i \ \ \text{and} \ \ \nabla_{Th (q,p)} (z) = \sum_{i = 0}^\infty b_i z^i $$ where $$\begin{cases} 
      a_i = b_i , & \text{for} \ i=0,1 \ (\text{mod} \ 4),\\
      a_i = -b_i , & \text{for} \ i=2,3 \ (\text{mod} \ 4).  
\end{cases}$$

As an immediate consequence of Problem~\ref{prob:polynomial}, one can compute the determinant of $Th (p,q)$ in general, which would be interesting to compare it with the generalized Kauffman--Harary conjecture in {\sc\S}\ref{sec:determinants}. Notice also that Problem~\ref{prob:polynomial} and Problem~\ref{prob:homological} are equivalent to each other due to the work of E.-S. Lee \cite{Lee05} and Ozsv\'ath and Szab\'o \cite{OS08_alexander}.

\begin{introprob}
\label{prob:homological}
Compute the following invariants of $Th (p,q)$ in terms of $p$ and $q$ more explicitly:
\begin{itemize}
    \item link Floer homology,
    \item Khovanov homology.
\end{itemize}    
\end{introprob}

Recall that Fox’s trapezoidal conjecture asserts that given an alternating knot $K$ with the Alexander polynomial $\Delta_K(t)=\pm \displaystyle\sum_{i=0}^{2n}\alpha_i(-t)^i$ with $\alpha_i>0$, there exists an integer $r \leq n$ such that: $$\alpha_0<\alpha_1< \dots < \alpha_{n-r}=\dots =\alpha_{n+r}> \dots ...>\alpha_{2n-1}>\alpha_{2n}.$$ We know that the conjecture is true for the Turk's head knot $Th (3,q)$, see {\sc\S}\ref{sec:alexander-conway}. The following problem is about the general case. One can consult Chbili's survey \cite{Chb24} for the current state of the art.

\begin{introprob}
Verify Fox’s trapezoidal conjecture for all Turk's head knots $Th (p,q)$.    
\end{introprob}

For alternating knots, recall that \emph{Hoste's condition} posits that the zeroes of the Alexander polynomials satisfy the equation: $$\mathfrak{Re}(z) > -1 .$$ In {\sc\S}\ref{sec:alexander-conway}, we saw that Hoste's condition holds for the Turk's head knots $Th (3,q)$. More generally, we ask the following open problem.

\begin{introprob}
Verify Hoste's condition for all Turk's head knots $Th (p,q)$.    
\end{introprob}

\subsection{On Hyperbolic Properties}

The \emph{volume}/\emph{determinant conjecture} (vol/det for short) posits that for any alternating hyperbolic link $L$ we have $$ \mathrm{vol} (L) < 2\pi \log \mathrm{det} (L) .$$ This conjecture was posed by Champanerkar, Kofman, and Purcell \cite{CKP16a} and it has been verified for the Turk's head links $Th(p,q)$ for $3 \leq p \leq 50$ and $2 \leq q \leq 50$. Combining the inequality (\ref{ineq:fibonacci}) in {\sc\S}\ref{sec:hyperbolic} with the equality (\ref{eq:determinant}) in {\sc\S}\ref{sec:determinants2}, we conclude that the vol/det conjecture holds for $Th(3,q)$. We wonder about all the remaining cases and ask: 

\begin{introprob}
\label{prob:vol/det}
Let $p \geq 4$. Verify the vol/det conjecture for all Turk's head links $Th (p,q)$.
\end{introprob}

Following Adams, Kaplan-Kelly, Moore, Shapiro, Sridhar, and Wakefield \cite{AKMSSW18}, the \emph{cusp crossing density} $d_{cc} (L)$ of a hyperbolic link is defined as $$d_{cc} (L) = \frac{cv (L)}{c (L)} \CommaPunct$$ where $cv (L)$ is the \emph{the maximal cusp volume} of $L$, i.e., it is the sum of the volumes of the individual cusps that do not overlap in their interiors and give the largest possible total.

In \cite{AKMSSW18}, they made various experiments by using SnapPy \cite{snappy} and raised the following interesting problem. They expect that for any link there is a Turk's head link that has a higher cusp crossing density.

\begin{introprob}
Verify whether the Turk's head links have the highest cusp crossing density of all the links.
\end{introprob}

More recently, Paoluzzi \cite{Pao21} proved a powerful result showing that prime alternating knots in $S^3$ are determined by any of their $n$-fold cyclic branched covers provided that $n>2$. We ask the following complementary problem by addressing the Turk's head knots:

\begin{introprob}
Is the double branched cover of $S^3$ branched along a Turk's head knot $Th (p,q)$ always hyperbolic? Does the double branched cover determine $Th (p,q)$?
\end{introprob}

In {\sc\S}\ref{sec:hyperbolic}, we mentioned about interesting observations of Mishra and Staffeldt relating the ranks of Khovanov homology groups with some density functions. We curiously raise the following problem.

\begin{introprob}
\label{prob:MS}
Explore the potential relations between twisted numbers, Khovanov ranks and volumes of the Turk's head knots $Th (p,q)$ more precisely.
\end{introprob}

\subsection{On Sliceness and Concordance}
The knot $K$ is called \emph{ribbon} if $K$ bounds an immersed disk in $S^3$ with only ribbon singularities. These singularities can be removed by pushing the disk into $B^4$, providing a slice disk $D$ for $K$ built only with $2$-dimensional $0$-, $1$-handles with respect to the radial Morse function on $B^4$.
Every ribbon knot is clearly slice by definition. The \emph{slice-ribbon conjecture}, an outstanding open problem posed by Fox \cite{Fox62}, posits that every slice knot is ribbon. It is well-known that the slice-ribbon conjecture was confirmed for $2$-bridge knots by Lisca \cite{Lis07, Lis07b} and for most pretzel and Montesinos knots by Greene and Jabuka \cite{GJ11} and Lecuona \cite{Lec12, Lec15, Lec18, Lec19}. 

Conway and Trotter \cite{Con70} studied the slice-ribbon conjecture by addressing Turk's head knots $Th (p,q)$ for odd values of $p$ and $q$, and they emphasized that potential counterexamples may be found among Turk's head knots. Thus, we would like to pose the following problem. Since Turk's head knots constitute a fairly interesting family of knots, it would be great to test the limits of this conjecture.

\begin{introprob}
\label{prob:conway-trotter}
Verify the slice-ribbon conjecture for the Turk's head knots $Th (p,q)$ when $p$ and $q$ are odd and $\mathrm{gcd} (p,q) =1$.
\end{introprob}

Recall that if $p$ is odd and $\mathrm{gcd} (p,q) =1$ then $Th(p,q)$ is strongly negative amphichiral. Assuming further $q$ to be odd, we know that they are also strongly positive amphichiral. In this case, various concordance invariants do not provide obstructions for sliceness, see the discussion in {\sc\S}\ref{sec:sliceness} and references therein. 

In {\sc\S}\ref{sec:sliceness}, we broadly discussed the known results about the orders of Turk's head knots $Th (3,q)$ in the concordance group $\mathcal{C}$. For the remaining cases, we raise the following problem:

\begin{introprob}
\label{prob:sliceness}
Determine the orders of Turk's head knots $Th (3,q)$ in $\mathcal{C}$ for $q \in 2\Z + 1 \setminus \{5,7,11,17,23\}$.
\end{introprob}

In {\sc\S}\ref{sec:sliceness}, we also explained that the knots $Th (3,q)$ for even values of $q$ are of order $2$ in $\mathcal{C}$. Since they are $\Q$-slice, they all non-trivially lie in $\mathrm{Ker} ( \psi: \mathcal{C} \to \mathcal{C}_\Q )$. In the following problem, we ask for their linear independence in this kernel.

\begin{introprob}
Let $p$ be even. Determine whether the Turk's head knots $Th (3,q)$ are $\Z_2$ linearly independent in $\mathrm{Ker} ( \psi)$.
\end{introprob}

The existence of a $\Z_2^\infty$ subgroup in $\mathrm{Ker} ( \psi)$ was previously proved by Cha \cite[Theorem~4.14]{Cha07}. His argument could be mimicked by understanding the Alexander polynomials of $Th (3,q)$ and their roots, see {\sc\S}\ref{sec:alexander-conway}. This was partially done in \cite{DPS24}.

We ask the following general concordance classification problem regarding the finite order elements in the concordance group $\mathcal{C}$.

\begin{introprob}
\label{prob:sliceness2}
Let $p$ be odd. Determine the orders of $Th (p,q)$ in $\mathcal{C}$.
\end{introprob}

For both Problem~\ref{prob:sliceness} and Problem~\ref{prob:sliceness2}, a reasonable strategy would be to understand the behavior of the twisted Alexander polynomial obstructions under Murasugi sums. This strategy was elaborated by Kindred \cite{Kin18} and Cheng, Hedden and Sarkar \cite{CHS22} for Khovanov homology and knot Floer homology, respectively. However, we cannot use $\tau$- and $s$-invariants here, since they are both equal to $\sigma$ up to some constant for alternating knots.

When $p$ is even, we know that the Trotter--Murasugi signature $\sigma$ is non-zero for $Th (p,q)$, see {\sc\S}\ref{sec:signatures}. So for a fixed $p$, the knot $Th (p,q)$ has infinite order in $\C$. As the next level, we ask the following concordance classification problem on infinite order elements. In this case, Levine-Tristam signatures $\sigma_\omega$ \cite{Lev69a, Tri69} might be useful, see also Takemura's divisibility criterion in Theorem~\ref{thm:takemura}. 

\begin{introprob}
Let $p$ be even. Prove that all the Turk's head knots $Th (p,q)$ are linearly independent in $\mathcal{C}$.
\end{introprob}

Due to the work of Cha and Ko \cite{CK02}, we know that the Trotter--Murasugi signature $\sigma$ also provides obstruction for $\Q$-sliceness. Similarly, for a fixed even value of $p$, the knot $Th (p,q)$ has infinite order in $\mathcal{C}_\Q$. Focusing on infinite order elements, we similarly ask the following $\Q$-concordance classification problem. Again, one can try to use Levine-Tristam signatures $\sigma_\omega$ to do so, since they factor through the $\Q$-concordance group $\mathcal{C}_\Q$, see again \cite{CK02}.

\begin{introprob}
Let $p$ be even. Prove that all the knots $Th (p,q)$ are linearly independent in $\mathcal{C}_\Q$.
\end{introprob}

Recall that all the Turk's head knots $Th (p,q)$ are strongly negative amphichiral for $p$ odd. Therefore, they are all $\Q$-slice and have order $1$ in $\mathcal{C}_\Q$ by Kawauchi's results \cite{Kaw80, Kaw09}.

\subsection{On Unknotting Numbers and Slice Genera}

A crucial invariant for a knot $K$ is  the \emph{slice genus} or \emph{$4$-genus}, which is defined as the minimal genus of an oriented surface embedded in $B^4$ whose boundary is $K$. The slice genus of $K$ is denoted by $g_4 (K)$. The following inequalities are well-known; see, for instance \cite{Liv05, Tan09}: $$\frac{\vert \sigma (K) \vert}{2} \leq g_4 (K) \leq u(K), \ \ \ \ \ \ g_4 (K) \leq g_3 (K), \ \ \ \ \ \ u (K) \leq \frac{c(K) -1}{2} \cdot$$ There is no inequality between the Seifert genus and the unknotting number. In particular, there are knots with $g_3 < u$ or $u < g_3$, see KnotInfo \cite{knotinfo}.

Focusing on Turk's head knots, we have the following formulas, coming from the above inequalities: $$\begin{cases} 
      \frac{q-1}{2} \leq g_4 (K) \leq u(K), & \text{if} \ p \ \text{is even,}\\
      0 \leq g_4 (K) \leq u(K) , & \text{if} \ p \ \text{is odd,}  
\end{cases}  \ \ \ \ \ \ g_4 (K) \leq \frac{(p-1)(q-1)}{2} \CommaPunct \ \ \ \ \ \ u (K) \leq \frac{(p-1)q -1}{2} \cdot$$

We identify the Turk's head knots with small crossing number by using SnapPy \cite{snappy}. Relying on KnotInfo \cite{knotinfo}, we have the following table, comparing their slice genera, unknotting numbers and Seifert genera. Note that the knots $Th(3,5)$, $Th(5,3)$, $Th(7,2)$ are slice (see {\sc\S}\ref{sec:surgeries}), so they both clearly have $g_4 = 0$.

\begin{table}[htbp]
\centering
\renewcommand{\arraystretch}{1.1}
\begin{tabular}{|l|l|l|l|}
\hline

Knots/Invariants & $g_4$ & $u$ & $g_3$  \\ \hline

$Th (3,2) = 4_1 $  & $1$ & $1$ & $1$  \\ \hline

$Th (3,4) = 8_{18} $  & $1$ & $2$ & $3$  \\ \hline

$Th (3,5) = 10_{123} $  & $0$ & $2$ & $4$  \\ \hline

$Th (4,3) = 9_{40} $  & $1$ & $2$ & $3$  \\ \hline

$Th (5,2) = 8_{12} $  & $1$ & $2$ & $2$  \\ \hline

$Th (7,2) = 12_{a477} $  & $0$ & $[2,3]$ & $3$  \\ \hline

$Th (5,3) = 12_{a1019} $  & $0$ & $2$ & $4$  \\ \hline

\end{tabular}
\vspace{0.3 em}
\caption{A comparison for invariants of Turk's head knots.}
\label{tab:comparison}
\end{table}

The order we use in Table~\ref{tab:comparison} is not random, since we know that the inequalities $$ g_4 (Th (3,q) ) \leq u (Th (3,q) ) \leq g_3 (Th (3,q) ) $$ hold due to the article by E.-K. Lee and S.-J. Lee \cite{LL13}. Taking the risk of having a negative answer, we propose the following conjecture:

\begin{introconj}
For the Turk's head knot $Th (p,q)$, the inequality $$u (Th (p,q) ) \leq g_3 (Th (p,q) ) $$ holds, and $g_4(Th(p,q))=u(Th(p,q))=g_3(Th(p,q))$ only for $Th (3,2)$.
\end{introconj}

Since the slice genus provides a lower bound for the unknotting number,\footnote{Compare with the \emph{untwisting number} of a knot defined to be the minimum number of \emph{null-homologous twists} required to convert a knot to the unknot. Allen, İnce, S. Kim, Ruppik and Turner \cite{All24} recently proved that the untwisting number gives an upper bound on the \emph{topological} slice genus.} the resolution of the following problem gives information about the above conjecture.

\begin{introprob}
For the Turk's head knot $Th (p,q)$, determine  $g_4 (Th (p,q) )$.    
\end{introprob}

Unfortunately, we cannot propose a tractable formula for the unknotting numbers. We expect this problem to be very difficult, and it might not be possible to find concrete patterns in terms of $p$ and $q$. Since Turk's head knots are closely related to torus knots, one can try to use unknotting sequences of torus knots to find some bounds, see work of Baader \cite{Baa10} and Siwach and Madeti \cite{SM15}.

\begin{introprob}
\label{prob:unknotting}
Determine the unknotting numbers of $Th (p,q)$ in terms of $p$ and $q$.
\end{introprob}

The \emph{unlinking number}, $u(L)$, for a link $L$, is similarly defined as the minimum number of crossing changes required to convert the diagram of $L$ into a diagram of an unlink. Presumably, the following problem is more difficult than Problem~\ref{prob:unknotting}.

\begin{introprob}
\label{prob:unlinking}
Determine the unlinking numbers of $Th (p,q)$ in terms of $p$ and $q$.
\end{introprob}

For several interesting results about unlinking numbers, see the articles by Kohn \cite{Koh93}, Motegi \cite{Mot96}, Jablan and Sazdanovi{\'c} \cite{JS07}, and Nagel and Owens \cite{NO15}.

\subsection{On Surgeries and Fillings}

In {\sc\S}\ref{sec:surgeries}, we observed that infinitely many surgeries along Turk's head knots give rise to homology $3$-spheres with $\Q$-homology $4$-ball fillings. Besides, we gave a few examples of homology $3$-spheres bounding contractible $4$-manifolds. We wonder about the general classification. 

\begin{introprob}
Classify integer and rational surgeries of Turk's head knots $Th(p,q)$ that yield homology $3$-spheres bounding contractible $4$-manifolds and $\Q$-homology $4$-balls.
\end{introprob}

\subsection{On Legendrian Properties and Invariants}

Recall that $Th(3,2) = 4_1$. In \cite{EH01}, Etnyre and Honda also proved that Legendrian figure-eight knots are determined by their maximal Thurston--Bennequin and rotation numbers, see also \cite{Etn05}. We wonder whether such a classification result can be achieved for Turk's head links $Th(p,q)$ in general.

\begin{introprob}
Classify the Legendrian Turk's head knots and links. Are they completely determined by their link types, maximal Thurston--Bennequin numbers, and rotation numbers?
\end{introprob}

In {\sc\S}\ref{sec:legendrian}, we saw that Turk's head links cannot bound orientable exact Lagrangian surfaces in the standard symplectic $4$-ball $(B^4, \omega_{std})$. However, Chen, Crider-Phillips, Reinoso, Sabloff and Yao \cite{CCRSY24} proved that the figure-eight knot $Th(3,2)$ bounds a non-orientable exact Lagrangian surface in $(B^4, \omega_{std})$, namely a Lagrangian Klein bottle with normal Euler number $4$. Inspired by their result and the work of Casals and Gao \cite{CG22}, we raise the following optimistic problems:

\begin{introprob}
Determine which Turk's head link $Th(p,q)$ has a non-orientable exact Lagrangian filling in $(B^4, \omega_{std})$.
\end{introprob}

\begin{introprob}
Is there any Turk's head knot or link $Th(p,q)$ that has more than one non-orientable exact Lagrangian fillings in $(B^4, \omega_{std})$?
\end{introprob}

\subsection{On Lissajous Knots} A \emph{Lissajous knot} $K(t)$ is a knot in $\mathbb{R}^3$ defined by the three parametric equations $$K(t)=\left (\cos (n_x t+\phi_x), \cos (n_y t+\phi_y), \cos (n_z t+\phi_z)\right ), \quad \text{for} \ 0 \leq t \leq 2\pi ,$$ where the frequencies $n_x, n_y,$ and $n_z$ are pairwise relatively prime integers and the phase shifts $\phi_x, \phi_y$ and $\phi_z$ are real numbers.

In \cite{BHJS94}, Bogle, Hearst, Jones and Stoilov introduced Lissajous knots motivated by Vassiliev theory, the physics of harmonic oscillators, and the biology of DNA dynamics. It is well-known that Lissajous knots are highly symmetric; in particular, they are all strongly positive amphichiral when $n_x, n_y$ and $n_z$ are odd.

Jones and Przytycki \cite{JP98} conjectured that the converse of the above implication is false by considering Turk's head knots $Th(3,q)$ for any odd $q$. Later, this conjecture appeared in Kirby's 1997 problem list \cite[Problem~1.62]{Kir78}; it is still an open problem.

Under the assumption that $n_x, n_y$ and $n_z$ are all odd, there are two main constraints for being Lissajous knots. Their Alexander polynomials must be perfect squares by the work of Hartley and Kawauchi \cite{HK79}. Moreover, their Arf invariants must be zero \cite{BHJS94}. This observation also follows from the fact that they are algebraically slice \cite{Lon84}. See the articles by Lamm \cite{Lam97}, Boocher, Daigle, Hoste, and Zheng \cite{BDHZ09}, and Soret and Ville \cite{SV16} for further discussions.

Recall that Turk's head knots $Th(p,q)$ are all strongly positive amphichiral when $p$ and $q$ are odd, see {\sc\S}\ref{sec:symmetries}. Generalizing the conjecture of Jones and Przytycki, we finally ask:

\begin{introconj}
Let $p$ and $q$ be odd. Then the Turk's head knots $Th (p,q)$ are not Lissajous knots.
\end{introconj}


\section{Appendix on Torus Knots and Links}
\label{sec:appendix}

The torus links have been studied extensively and explored from various perspectives. We follow Murasugi's book \cite[{\sc\S}7]{Mur96} as the main reference, so we will only cite the other references as needed.

Let $p$ and $q$ be two positive integers. Recall from {\sc\S}\ref{sec:introduction} that the \emph{torus link} $T(p,q)$ is given by $$ T (p,q) \doteq \reallywidehat{ ( \sigma_1 \sigma_2 \sigma_3 \sigma_4  \ldots \sigma_{p-1} )^q } . $$ The torus link is also realized as the link of the Brauner type complex curve singularity $$\Sigma(p,q) \doteq \left \{ (x,y) \in \mathbb{C}^2 \ \vert \  x^p + y^q = 0 \right \} \cap S^3 \subset \mathbb{C}^2 ,$$ see Milnor's book and references therein \cite[{\sc\S}1]{Mil68}.

The torus link has $\mathrm{gcd}(p,q)$ components, hence it is a knot and said to be the \emph{torus knot} if $\mathrm{gcd}(p,q) = 1$. Because of the braid convention, the torus link $T(p,q)$ is called the \emph{positive torus link} or the \emph{right-handed torus link}. The link diagram for $T(p,q)$ depicted in Figure~\ref{fig:torusknots} is minimal, and the minimum number of crossings of $T(p,q)$ is given by $$ c(T(p,q)) = \min \{ p(q-1), q(p-1) \} = pq - \max \{p,q\} .$$

Clearly, $T(1,q)$ or $T(p,1)$ gives the unknot for any $p$ and $q$. So from now on we assume that $p \geq 2$ and $q \geq 2$. Unlike the Turk's head link $Th (p,q)$, the link types of $T(p,q)$ and $T(q,p)$ are the same, i.e., they are always isotopic to each other. The mirror image of $T(p,q)$ is given by $$ \overline{T (p,q)} \doteq \reallywidehat{ ( \sigma^{-1}_1 \sigma^{-1}_2 \sigma^{-1}_3 \sigma^{-1}_4  \ldots \sigma^{-1}_{p-1} )^q }. $$ In a similar vein, the link $ \overline{T(p,q)} $ is said to be the \emph{negative torus link} or the \emph{left-handed torus link}.

The torus link $T(p,q)$ is never amphichiral, i.e., the link types of $T(p,q)$ and $\overline{T (p,q)}$ are always different. However, the link $T(p,q)$ is invertible, i.e., $T(p,q)$ is isotopic to $-T(p,q)$. Moreover, $T(p,q)$ is both $p$-periodic and $q$-periodic. In particular, the torus link $T(p,q)$ is strongly invertible and $$\mathrm{Sym} (S^3, T(p,q) ) \cong \Z_2 ,$$ see \cite[{\sc\S}10.6]{Kaw90}.

It is well-known that the torus link $T(p,q)$ is always fibered. Thus, $T (p,q)$ bounds a unique minimal genus Seifert surface rel. boundary up to isotopy in $S^3$, see Kobayashi's work \cite{Kob89}. The Seifert genus of $T (p,q)$ is given by $$ g_3 ( T(p,q) ) = \frac{ (p-1)(q-1) + \mathrm{gcd} (p,q) - 1 }{2} \cdot $$ Moreover, the Seifert matrix $V_{pq}$ associated to torus link $T (p,q)$ is given by the following square matrix of size $(p-1)(q-1)$: $$ V_{p,q} = \begin{pNiceMatrix}
 A_q    &            &        &       &        &        \\
-A_q   & A_q     &        &        &        &        \\
       & -A_q       & A_q    &        &        &        \\
       &            & \ddots & \ddots &        &        \\
       &            &        & -A_q   & A_q &        \\
       &            &        &        & -A_q   & A_q    \\
\end{pNiceMatrix} \CommaPunct $$ where $A_q$ is the Seifert matrix for $T(2,q)$ of size $q-1$ (cf.~{\sc\S}\ref{sec:alexander-general}):

$$A_q = \begin{pNiceMatrix}
-1     & 1      &        &        &        \\
       & -1     & 1      &        &        \\
       &        & \ddots & \ddots &        \\
       &        &        & \ddots & 1      \\
       &        &        &        & -1     \\       
\end{pNiceMatrix}.$$ Here, all the other blocks are zero blocks of size $q-1$. Therefore, the Alexander polynomial of the torus link is given by $$ \Delta_{T(p,q)} (t) = \frac{(1-t)(1-t^\frac{pq} {d})^d}{(1-t^p)(1-t^q)} \CommaPunct $$ where $d = \mathrm{gcd} (p,q)$. Due to the articles by Jones \cite{Jon87} and Isodro, Labastida and Ramallo \cite{ILR93}, the Jones polynomial of $T(p,q)$ is completely determined:

$$V_{Th (p,q)} (s) = \frac{1}{1-s^2} \sum_{i=0}^{d} \begin{pmatrix}
			d \\
			i
		\end{pmatrix} s^{ \frac{p}{d} \left ( 1 + \frac{q}{d} \right ) (d-i)} \left ( s^{ \frac{q}{d} (d-i) } - s^{ 1 + \frac{q}{d} i } \right ) \cdot$$

\noindent The determinant of the torus link $T(p,q)$ can be computed in terms of $p$ and $q$ from either its Alexander or Jones polynomial. In particular, for the torus knots we have $$ \mathrm{det} (T(p,q)) = \begin{cases} 
      1, & \text{if} \ p \ \text{and} \ q \ \text{are both odd,} \\
      p, & \text{if} \ p \ \text{is odd and} \ q \ \text{is even.}
\end{cases}$$ Given a prime number $n$, $T(p,q)$ is therefore not $n$-colorable if $p$ and $q$ are both odd. On the other hand, if $p$ is odd, $q$ is even and $n \vert p$, then $T(p,q)$ is $n$-colorable. Moreover, Breiland, Oesper and Taalman \cite{BOT09} proved that every $n$-colorable torus knot has exactly one nontrivial $n$-coloring class. Furthermore, the lower and upper bounds for the minimal number of colorings of torus knots and links were studied by Kauffman and Lopes \cite{KL08}, Jablan, Kauffman and Lopes \cite{JKL13}, and Ichihara and Matsudo \cite{IM17}.

Comparing the ones for Turk's head links $Th(p,q)$ that appeared in {\sc\S}\ref{sec:signatures}, one can obtain a Murasugi sum decomposition of torus links $T(p,q)$ in the following fashion: $$T(p,q) = \underbrace{T(2,q) \ \#_{2q} \ T(2,q) \ \#_{2q} \ \ldots \ \#_{2q} \  T(2,q) }_{(p-1) \ \text{copies}} . $$

Gordon, Litherland and Murasugi \cite{GLM81} provided the following recursive formula for the Trotter--Murasugi signatures of torus links: \begin{itemize}
    \item[I.] Let $2q < p$. Then \begin{itemize}
        \item[$\bullet$] if $q$ is odd, then $\sigma ( T(p,q) ) - \sigma ( T(p-2q,q) ) = -(q^2 -1)$,
        \item[$\bullet$] if $q$ is even, then $\sigma ( T(p,q) ) - \sigma ( T(p-2q,q) ) = -q^2$,
    \end{itemize}

    \item[II.] Let $q \leq p < 2q$. Then \begin{itemize}
        \item[$\bullet$] if $q$ is odd, then $\sigma ( T(p,q) ) + \sigma ( T(2q-p,q) ) = -(q^2 -1)$,
        \item[$\bullet$] if $q$ is even, then $\sigma ( T(p,q) ) + \sigma ( T(2q-p,q) ) = -(q^2 -2)$,
    \end{itemize}

    \item[III.] $\sigma ( T(2q,q) ) = -(q^2 -1) $,

    \item[IV.] $\sigma ( T(p,q) )= \sigma ( T(q,p))$, $\sigma ( T(1,q)) = 0$ and $\sigma ( T(2,q)) = -(q-1)$.
\end{itemize}

\noindent Using Levine-Tristam signatures, Litherland \cite{Lit79} proved that the family of all torus knots $T(p,q)$ are linearly independent in the concordance group $\mathcal{C}$. Furthermore, they generate a $\Z^\infty$ summand in $\mathcal{C}$.

The work of Kroinheimer and Mrowka \cite{KM93, KM95} resolved the Milnor conjecture affirmatively and proved the relationship between the Seifert genera, the slice genera, and the unlinking numbers for all the torus links: $$ g_3 ( T(p,q) ) = g_4 ( T(p,q) ) = u ( T(p,q) ) = \frac{ (p-1)(q-1) + \mathrm{gcd} (p,q) - 1 }{2} \cdot $$

As Legendrian links, torus knots and links were well-studied. In \cite{EH01}, Etnyre and Honda computed their maximal Thurston--Bennequin numbers: $$ \TB (T(p,q)) = pq - p - q \quad \text{and} \quad \TB (\overline{T(p,q)}) = -pq . $$ Moreover, they classified the Legendrian torus knots by using the technique of convex surfaces. More precisely, Etnyre and Honda proved that all positive and negative torus knots are determined up to Legendrian isotopy by their knot types, Thurston--Bennequin invariants, and rotation numbers. The recent work of Dalton, Etnyre and Traynor \cite{DET24} achieved the classification of all Legendrian torus links.

Another important work for the study of Legendrian torus links is the article by Casals and Gao \cite{CG22} in which they proved that every maximal Thurston--Bennequin number positive Legendrian torus link $T(p,q)$ admits infinitely many exact Lagrangian fillings in the standard symplectic $4$-ball $(B^4, \omega_{std})$, except for $T(2,q)$, $T(3,3)$, $T(3,4)$, and $T(3,5)$.

Now, we discuss Dehn surgeries along torus knots, which are always Seifert-fibered spaces with three singular fibers. This classical result was due to Moser \cite{Mos71}, but we follow the standard convention appearing in the article by Owens and Strle \cite[Lemma~4.4]{OS12}. For any rational number $r$, we have $$S^3_r(T(p,q))=Y\left(2;\frac{p}{q^*} \CommaPunct \frac{q}{p^*} \CommaPunct \frac{pq-r}{pq-r-1}\right) \CommaPunct$$ where $q^*$ is the multiplicative inverse of $q$ modulo $p$, i.e. $qq^*\equiv 1\pmod p$ and $1\leq q^*<p$. Similarly, for $p$, we exchange the roles of $p$ and $q$. In particular, we have the following identifications: \begin{itemize}
    \item $S^3_{pq \pm 1 }(T(p,q)) = L (pq \pm 1, \mp q^2) = L(pq \pm 1, \mp p^2)$,
    \item $S^3_{pq }(T(p,q)) = L (p,-q) \# L (q,-p)$,
    \item $S^3_{-1/n }(T(p,q)) = \Sigma(p,q,pqn+1)$ for $n \geq 1$,
    \item $S^3_{-1/n }(\overline{T(p,q)}) = \Sigma(p,q,pqn-1)$ for $n \geq 1$,
\end{itemize}

\noindent where $L(a,b)$ and $\Sigma(a,b,c)$ respectively denote \emph{lens spaces} and \emph{Brieskorn homology $3$-spheres}, which are defined and oriented as follows: $$L(a,b) \doteq S^3 / \Z_a, \quad \text{where} \ S^3 \subset \mathbb{C}^2 \ \text{and} \ (z_1, z_2) \mapsto (e^{2\pi i / a}, e^{2\pi ib / a}), $$ and $$\Sigma(a,b,c) \doteq \left \{ (x,y,z) \in \mathbb{C}^3 \ \vert \  x^a + y^b + z^c = 0 \right \} \cap S^5 \subset \mathbb{C}^3 .$$ 

Finally, we mention that Aceto, Golla, Larson, and Lecuona classified all positive integer surgeries of torus knots $T(p,q)$ bounding $\Q$-homology $4$-balls. See their article \cite{AGLL20} for the complete classification and further discussion.


\section*{Afterword}

We started typing this survey during the course of writing our article \cite{DPS24}. Using KnotInfo \cite{knotinfo} and SnapPy \cite{snappy}, the Turk's head knots $Th (3,q)$ first attracted the authors because of their rich symmetries. 

After tracing several reference articles, we realized that these types of knots and links appear in the literature under different names: Turk's head knots/links, rosette knots/links and weaving knots/links. Along the way, we also observed that several results on $Th (p,q)$ are either proved or observed simultaneously by different authors, so we decided to write this survey to share common knowledge with the community. This allowed us to review several important classical articles in the field of knot theory and low-dimensional topology as well. 

We finally noted that there are several families of knots and links in the literature that are direct or indirect generalizations of Turk's head knots and links, such as \emph{spiral knots and links} \cite{BETV10, KST16}, \emph{hybrid or generalized weaving links} \cite{SMR21, AC23}, and \emph{Lissajous toric knots} \cite{SV20}. A curious reader may consult these references for their definitions and interesting properties.


\section*{Acknowledgements}

We greatly benefited from the constructive and valuable feedback provided by M{\'a}rton Beke, Carlo Collari, Vikt\'oria F\"oldv\'ari, Marco Golla, Christoph Lamm, Brendan Owens, J{\'o}zef H. Przytycki, and Andrei Yu. Vesnin. We would like to thank Keegan Boyle, Ana G. Lecuona, Paolo Lisca, Charles Livingston, Makoto Sakuma, and Andr\'as I. Stipsicz for their encouraging comments. We are grateful to Charles Livingston, Luisa Paoluzzi and Daniel Ruberman for their helpful comments regarding the existing literature. We would also like to thank Darren D. Long and Yasutaka Nakanishi for sharing the digital copies of \cite{Lon83} and \cite{NY00}, and Andr\'as I. Stipsicz for providing us with the scans of \cite{Boz80, Boz85}.

The first author thanks the Italian National Group for Algebraic and Geometric Structures and their Applications (INdAM-GNSAGA). The second author was supported by the CNRS postdoctoral fellowship at the Laboratoire de Math\'ematiques Jean Leray in Nantes Universit\'e, France.

Finally, we would like to thank the anonymous referee for their careful assessment and invaluable feedback.

\bibliography{references}
\bibliographystyle{amsalpha}

\end{document}